\documentclass[11pt]{amsart}
\usepackage[utf8]{inputenc}
\usepackage[english]{babel}

\usepackage{amscd,amssymb}
\usepackage{amsmath}
\usepackage{amsthm}
\usepackage{graphicx}
\usepackage{fancybox}
\usepackage{epic,eepic}
\usepackage{amstext}
\usepackage{mathtools}
\usepackage{esint}
\usepackage[square,numbers]{natbib}
\usepackage{booktabs}
\usepackage[hide links]{hyperref}
\usepackage{comment}
\usepackage{mathtools}
\usepackage{tikz,pgfplots}
\usepackage{subcaption}

\usepackage[noabbrev, capitalise, nameinlink]{cleveref}
\crefname{equation}{}{}
\crefname{assumption}{Assumption}{Assumptions}
\crefformat{equation}{\textup{#2(#1)#3}}

\usepackage{algorithm}
\usepackage{algpseudocode}
\usepackage{tabularx}
\newcolumntype{Y}{>{\centering\arraybackslash}X}

\newtheorem{theorem}{Theorem}[section]

\newtheorem{lemma}[theorem]{Lemma}
\newtheorem{assumption}[theorem]{Assumption}

\newtheorem{conjecture}[theorem]{Conjecture}
\theoremstyle{definition}

\theoremstyle{remark}
\newtheorem{remark}[theorem]{Remark}
\numberwithin{theorem}{section}
\numberwithin{equation}{section}
\numberwithin{table}{section}
\numberwithin{figure}{section}

\def\elem{T}

\def\B{\mathcal{B}}

\def\V{\mathcal{V}}
\def\R{\mathcal{R}}

\def\Pnull{\mathbb{P}_0}
\newcommand{\N}{\mathbb{N}}
\renewcommand{\bullet}{\cdot}

\def\T{\mathcal{T}}
\def\TH{\mathcal{T}_H}

\def\pset{\mathcal M}
\def\ptr{\mathcal M^\mathrm{tr}}
\def\prb{\mathcal M^\mathrm{rb}}

\newcommand{\with}{\,:\,}

\def\PiH{\Pi_H}

\def\linh{\mathrm{span}}

\newcommand{\A}{\mathcal{A}}

\def\one{\mathbf{1}}

\newcommand{\tnormf}[2]{\| #1 \|_{L^2(#2)}}

\newcommand{\vnormf}[2]{\| #1 \|_{\mathcal V,#2}}
\newcommand{\vnormof}[1]{\| #1 \|_{\mathcal V}}
\newcommand{\vspf}[3]{(#1,#2)_{\mathcal V,#3}}
\newcommand{\vspof}[2]{(#1,#2)_{\mathcal V}}

\newcommand{\tspf}[3]{{( #1,#2 )}_{L^2(#3)}}

\newcommand{\patch}{\mathsf{N}}

\DeclareMathOperator*{\argmin}{arg\,min}

\allowdisplaybreaks

\begin{document}

\begin{abstract}
	This paper presents a method for the numerical treatment of reaction-convection-diffusion problems with parameter-dependent coefficients that are arbitrary rough and possibly varying at a very fine scale. The presented technique combines the reduced basis (RB) framework with the recently proposed super-localized orthogonal decomposition (SLOD). More specifically, the RB is used for accelerating the typically costly SLOD basis computation, while the SLOD is employed for an efficient compression of the problem's solution operator requiring coarse solves only.
	The combined advantages of both methods allow one to tackle the challenges arising from parametric heterogeneous coefficients. Given a value of the parameter vector, the method outputs a corresponding compressed solution operator which can be used to efficiently treat multiple, possibly non-affine, right-hand sides at the same time, requiring only one coarse solve per right-hand side. 
\end{abstract}
	
\title[A reduced basis super-localized orthogonal decomposition]{A reduced basis super-localized orthogonal decomposition for reaction-convection-diffusion problems}
\author[F.~Bonizzoni, M.~Hauck, D.~Peterseim]{Francesca Bonizzoni$^*$, Moritz Hauck$^\dagger$, Daniel Peterseim$^{\dagger,\ddagger}$}
\address{${}^*$ MOX - Department of Mathematics, Politecnico di Milano, via Bonardi 9, 20133 Milano, Italy}
\email{francesca.bonizzoni@polimi.it}
\address{${}^{\dagger}$ Institute of Mathematics, University of Augsburg, Universit\"atsstr.~12a, 86159 Augsburg, Germany}
\email{moritz.hauck@uni-a.de}
\address{${}^{\ddagger}$ Institute of Mathematics \& Centre for Advanced Analytics and Predictive Sciences (CAAPS), University of Augsburg, Universit\"atsstr.~12a, 86159 Augsburg, Germany}
\email{daniel.peterseim@uni-a.de}
\thanks{The work of the authors is part of a project that has received funding from the European Research Council ERC under the European Union's Horizon 2020 research and innovation program (Grant agreement No.~865751). F. Bonizzoni is member of the INdAM Research group GNCS}

\maketitle

{\tiny {\bf Keywords.}~Parameter-dependent PDE, model order reduction, reduced basis, multiscale method,\\
	\hphantom{\indent {\tiny {\bf Keywords.}}}~numerical homogenization, reaction-convection-diffusion problems
}\\
\indent
{\tiny {\bf AMS subject classifications.}  
	{\bf 35B30}, 
	{\bf 41A65}, 
	{\bf 65N12},
	{\bf 65N15},
	{\bf 65N30}
} 
%

\section{Introduction}

Many phenomena in engineering and sciences can be modeled by partial differential equations (PDEs) with highly heterogeneous and strongly varying coefficients with oscillations appearing on several non-separated scales. In the literature, such PDEs are referred to as multiscale PDEs. 
Often, mathematical models also depend on a set of parameters accounting, e.g., for varying material properties or geometric variability, leading to parametrized multiscale PDEs. Typical examples include subsurface fluid flows in porous media and composite materials.
For this class of problems - even for fixed values of the parameter vector - the solution by means of classical finite element methods is challenging. Indeed, a very fine mesh is required to resolve all fine-scale features, resulting in large linear systems and possibly days of CPU time. If one is interested in the PDE solution corresponding to various parameter input values, the repeated numerical solution of the PDE may easily become computationally unfeasible.

Non-parametric multiscale PDEs have been intensively studied by the numerical analysis community over the last decades. For such problems, several kinds of methods have been established which we summarize under the term numerical homogenization. Under minimal assumptions on the coefficients, these methods are able to achieve uniform orders of convergence. This is achieved by using problem-adapted ansatz spaces possessing local, efficiently computable basis functions. However, there is always a computational overhead, i.e., either the support of the basis functions or the number of basis functions needs to be increased with the desired accuracy. Prominent methods achieving almost optimal numerical homogenization include the multiscale spectral generalized finite element method~\cite{BaL11,EFENDIEV2013116,Ma22}, the adaptive local basis~\cite{GGS12,Wey17}, the localized orthogonal decomposition method (LOD)~\cite{HeP13,MaP14,KPY18}, rough polyharmonic splines~\cite{OZB14}, and gamblets~\cite{Owh17}. 
For an overview on numerical homogenization, we refer to the textbooks~\cite{MalP20,OwhS19} and the recent review article~\cite{Peterseim2021}. 
Recently, the super-localized orthogonal decomposition (SLOD) has been proposed~\cite{HaPe21b} which is a version of the LOD with significantly improved localization properties (see also~\cite{Freese-Hauck-Peterseim,BFP22,graphSLOD,pumslod}).
Using its super-localized problem-adapted basis functions, the SLOD  achieves a particularly efficient compression of the solution operator of the PDE only requiring the solution of a sparse coarse linear system. For the construction of the SLOD basis functions, local PDEs need to be solved on the fine scale. While this step can be performed at moderate computational costs for non-parameteric PDEs, it may get unfeasible for parametric PDEs.

Thus, we handle the parameter dependence by adopting the model order reduction (MOR) perspective. Namely, the high-fidelity problem - also referred to as full order model (FOM) - is replaced by a suitable reduced order model (ROM), which is accurate and at the same time fast to evaluate (in particular, its cost is independent of the dimension of the original FOM). ROM techniques are two-phase approaches. First an offline training phase is performed, consisting in the numerical computation of FOM solutions corresponding to an appropriate set of input values of the parameter vector. The resulting set of solutions (also known as snapshots) are then employed for constructing an approximation of the parameter-to-solution map (the surrogate), which is then evaluated online at any new parameter values of interest. 
Between the various ROMs, the reduced basis (RB) method~\cite{QMN15,HesthavenRozzaStamm-2016} represents a remarkable technique to perform model reduction. It constructs the surrogate by projection over a (small) finite dimensional subset being the output of a principal component analysis or of a greedy search procedure.

The efficient treatment of parametric multiscale PDEs - as we are interested in here - calls for technologies that merge the features of model order reduction and numerical homogenization. 
During the last decades, there have been many works in this direction including the localized RB multiscale method~\cite{KaulmannOhlbergerHaasdonk-2011,KaulmannFlemischHaasdonkLieOhlberger-2015,OhlbergerSchindler-2014,OhlbergerSchindler-2015,AHK12,BIO21}, the RB-LOD which is a combination of the RB and LOD~\cite{AbH15,KeilRave-2021}, and ArbiLoMod~\cite{BuhrEngwerOhlbergerRave-2017}. Additionally, we mention the RB multiscale finite element method~\cite{Nguyen-08,HesthavenZhangZhu-2015} and the FE$^2$-based MOR method~\cite{HeAveryFarhat-2020} which, however, require scale separation of the coefficients. The key idea of the above-mentioned methods is to localize the snapshot computation, i.e., fine-scale snapshots are only computed locally on subdomains. This makes the snapshot computation feasible also for parametric multiscale PDEs. 

In this paper, we introduce the reduced basis super-localized orthogonal decomposition (RB-SLOD), a novel methodology for the numerical treatment of reaction-convection-diffusion problems with arbitrary rough and parameter-dependent coefficients. It is the result of the careful combination of the SLOD method - used for the localization of the PDE - with the RB method - used for the acceleration of local computations. The significantly improved localization properties of the SLOD allow to compute the snapshots on much smaller subdomains compared to the RB-LOD, the target accuracy being prescribed. This is the key to significant computational savings in both the offline and online phases. 
Given a value of the parameter vector, the proposed method outputs a corresponding compressed solution operator.
A key strength is that the compressed solution operator can be used for the efficient computation for various right-hand sides only requiring one coarse solve per right-hand side. In particular, the proposed algorithm is not affected by right-hand sides with a non-affine parameter dependence, while, typically, for such right-hand sides, methods like the empirical interpolation method~\cite{BMN04} need to be applied causing extra computational costs.
Moreover, it enables fast and accurate computations for any parameter value, making the algorithm suitable in multi-queries and real-time scenarios.

The performed numerical tests confirm the proposed method's effectiveness and highlight its advantages compared to other state-of-the-art methodologies. For a pure diffusion problem with a highly oscillating parametric diffusion coefficient, a comparison with the RB-LOD shows that, for reaching the same accuracy, the RB-SLOD requires much smaller oversampling domains. This results in smaller local patch problems and a sparser coarse system. Moreover, the effectiveness of the RB-SLOD is also demonstrated for a parametric mass transfer problem (i.e., a reaction-convection-diffusion problem) with a non-affine right-hand side. 

The paper is organized as follows. In~\cref{sec:problem_setting}, we formulate the problem of interest and introduce some notation. For a fixed parameter value,~\cref{sec:SLOD} then recalls the definition of the SLOD. The core of this paper are~\cref{sec:RB-SLOD,sec:error_analysis}, where the RB-SLOD method is introduced and analyzed. In a series of numerical experiments,~\cref{sec:numexp} underlines the method's effectiveness for parametric multiscale PDEs. Finally, conclusions are drawn in~\cref{sec:conclusions}.

\section{Problem setting}
\label{sec:problem_setting}

Let $\Omega \subset \mathbb{R}^d$, $d\in\{2,3\}$ be a polygonal/polyhedral Lipschitz domain which is scaled such that its diameter is of order one. Denoting by $\pset \subset \mathbb R^p$, $p \in \N$, a compact parameter set, we consider the parameter-dependent reaction-convection-diffusion problem posed in the domain $\Omega$ with parameters $\mu$ from $\pset$, i.e.,  
\begin{equation}
\label{eq:strongform}
-\operatorname{div}(a_\mu\nabla u_\mu) + b_\mu \cdot \nabla u_\mu + c_\mu u_\mu = f,
\end{equation}
where $c_\mu$, $b_\mu$, and $a_\mu$ denote the reaction, convection, and diffusion coefficients, respectively and $f\in L^2(\Omega)$ is the parameter-independent right-hand side.
We assume that
$c_\mu \in L^\infty(\Omega,\mathbb R)$ and $b_\mu \in L^\infty(\Omega,\mathbb R^d)$ with $L^\infty$-norms bounded uniformly in the parameter $\mu$. Further, let the matrix valued diffusion coefficient $a_\mu \in L^\infty(\Omega, \mathbb{R}^{d\times d})$ be symmetric and uniformly positive definite, i.e., there exist parameter-independent constants $0<\lambda\leq \Lambda < \infty$ such that, for all $\mu \in \pset$, $\eta \in \mathbb R^d$, and almost every $x \in \Omega$, it holds
\begin{equation}
\label{eq:unifboundsA}
\lambda |\eta|^2\leq \eta \cdot (a_\mu(x)\eta) \leq \Lambda |\eta|^2,
\end{equation}
with $|\cdot|$ denoting the Euclidean norm of $\mathbb R^d$. 

Given the pairwise disjoint partition of the boundary $\partial \Omega = \Gamma_1\cup \Gamma_2\cup \Gamma_3$ (with $\Gamma_3$ closed), we supplement~\cref{eq:strongform} with the following homogeneous mixed boundary conditions
\begin{equation}
\begin{aligned}
(a_\mu\nabla u_\mu) \cdot \nu &= 0 \quad &&\text{on }\Gamma_1,\\
(a_\mu\nabla u_\mu) \cdot \nu  + d_\mu u&= 0 \quad &&\text{on }\Gamma_2,\\
u_\mu&= 0\quad &&\text{on }\Gamma_3.
\end{aligned}
\label{eq:boundarycond}
\end{equation}
In~\cref{eq:boundarycond}, $\nu$ denotes the outer unit normal vector and $d_\mu \in L^\infty(\Gamma_2,\mathbb R)$ is  non-negative and its $L^\infty$-norm is uniformly bounded in~$\mu$. 

Let $\V \coloneqq \{v \in H^1(\Omega,\mathbb R) \with v|_{\Gamma_3} = 0\}$, where $v|_{\Gamma_3}$ shall be understood in the sense of traces. We equip the space $\V$ with the standard $H^1$-inner product and the corresponding induced norm, namely:
\begin{equation*}
\vspof{u}{v} \coloneqq \tspf{\nabla u}{\nabla v}{\Omega} + \tspf{u}{v}{\Omega},\quad \vnormof{u}^2 \coloneqq \vspof{u}{u}.
\end{equation*}
For any fixed parameter value $\mu \in \pset$, the weak formulation of equation~\cref{eq:strongform} with boundary conditions~\cref{eq:boundarycond} reads: find $u_\mu \in \V$ such that, for all $v \in \V$,
\begin{equation}
\label{eq:weak}
A_\mu(u_\mu,v) = \tspf{f}{v}{\Omega},
\end{equation} 
where the bilinear form $A_\mu\colon \V\times\V\to\mathbb R$ is given by
\begin{equation*}
A_\mu(u,v) \coloneqq \tspf{a_\mu \nabla u}{\nabla v}{\Omega} + \tspf{b_\mu\cdot \nabla u}{v}{\Omega} + \tspf{c_\mu u}{v}{\Omega} + \tspf{d_\mu u}{v}{\Gamma_2}.
\end{equation*}
Under standard assumption on the coefficients $d_\mu,c_\mu,$ and $b_\mu$ (see, e.g., the textbook~\cite[Ch.~3.2]{knabner2021numerical}) one can prove the coercivity and continuity of the bilinear form $A_\mu$
\begin{equation}
\label{eq:conditions4laxmilgram}
A_\mu(u,u)\geq \alpha \vnormof{u}^2,\quad A_\mu(u,v) \leq \beta \vnormof{u}\vnormof{v}.
\end{equation}
If conditions~\cref{eq:conditions4laxmilgram} are fulfilled, the well-posedness of the weak formulation~\cref{eq:weak} follows by standard arguments and the Lax--Milgram theorem.

As customary in the reduced basis context, we make the following assumption.
\begin{assumption}[Affine decomposition]
	The bilinear form~$A_\mu$ can be decomposed into~$Q$ terms as follows
	\begin{equation}
	\label{eq:decA}
	A_\mu(u,v) = \sum_{q=1}^{Q} \theta_q(\mu) B_q(u,v)
	\end{equation}
	with parameter-independent continuous bilinear forms $B_q\colon \V\times \V\to \mathbb R$ and measurable functions $\theta_q\colon \pset \to \mathbb R$.
\end{assumption}

\begin{remark}[Practical limitations]
	From the theoretical point of view, the RB-SLOD can be applied as long as the bilinear form $A_\mu$ fulfills~\cref{eq:decA}, i.e., independently of the actual value of $Q$. However, as usual for RB methods, there are limitations on the number of summands $Q$ in~\cref{eq:decA} and smoothness requirements with respect to $\mu$.
	Indeed, for large numbers $Q$ or lacking smoothness, the offline phase of the RB-SLOD gets increasingly expensive since a lot of precomputations are required. In such cases, the proposed RB-SLOD might not pay off, see also \cref{rem:decdelta}.
\end{remark}

\begin{remark}[Approximate affine decomposition]
	\label{rem:appaffdec}
	If no affine decomposition of $A_\mu$ in the sense of~\cref{eq:decA} is available, one can use the empirical interpolation method~\cite{BMN04} to compute affine approximations of the PDE coefficients $d_\mu,c_\mu,b_\mu,$ and $a_\mu$. The RB-SLOD can then be applied to these approximate coefficients. Note that the number of terms in the affine approximations of the coefficients depends on their smoothness with respect to $\mu$, i.e., for a prescribed accuracy, smoother coefficients can be approximated by sums with a smaller number of terms.
\end{remark}


For the ease of presentation, we introduce notation which is frequently used in the following sections. We denote the associated solution operator that maps a given right-hand side $f\in L^2(\Omega)$ to the corresponding unique solution $u_\mu$ by $\A_\mu^{-1}\colon L^2(\Omega)\to \V$. 
Moreover, for subsets $\omega\subset \Omega$, we denote by $A_{\omega,\mu}(\bullet,\bullet)$,  $B_{\omega,q}(\bullet,\bullet)$, $\vspf{\bullet}{\bullet}{\omega}$, and $\vnormf{\bullet}{\omega}$ the restrictions of $A_\mu(\bullet,\bullet)$, $B_{q}(\bullet,\bullet)$, $\vspof{\bullet}{\bullet}$, and $\vnormof{\bullet}$ to $\omega$, respectively. 
We denote the local solution space by
\begin{equation}
\label{eq:locsolspace}
\V_\omega \coloneqq \{v \in H^1(\omega)\with v|_{\partial \omega \backslash(\Gamma_1\cup \Gamma_2)} = 0\}
\end{equation}
and define the local solution operator $\mathcal A_{\omega,\mu}^{-1}\colon L^2(\omega)\to \V_\omega$ as the operator mapping any right-hand side $f_\omega \in L^2(\omega)$ to the local solution $u_{\omega,\mu}\in \V_\omega$ satisfying, for all $v \in \V_\omega$,
\begin{equation}
\label{eq:local_pde}
A_{\omega,\mu}(u_{\omega,\mu},v) = \tspf{f_\omega}{v}{\omega}.
\end{equation}

\section{SLOD for non-parametric reaction-convection-diffusion problems}
\label{sec:SLOD}

This section introduces the SLOD technique for non-parametric reaction-convection-diffusion problems with arbitrary rough PDE coefficients $d_\mu, c_\mu,b_\mu,$ and $a_\mu$. For this purpose, we fix the parameter $\mu \in \pset$ throughout this section.
We consider the possibly coarse, quasi-uniform triangular or quadrilateral mesh $\TH$ of $\Omega$ consisting of closed, convex, and shape-regular elements, where the subscript $H$ denotes the maximal element diameter. In general, the mesh $\TH$ does not resolve the small scale variations of the PDE coefficients.
Let $\Pnull(\TH)$ denote the space of $\TH$-piecewise constants 
\begin{equation*}
\Pnull(\TH)=\operatorname{span}\left\{\one_K\with K \in \TH\right\}
\end{equation*}
with $\one_K$ denoting the characteristic function of the element $K \in \TH$. Furthermore, by $\Pi_H\colon L^2(\Omega)\rightarrow \Pnull(\TH)$, we denote the $L^2$-orthogonal projection onto $\Pnull(\TH)$, namely, for all $K\in\TH$, $\Pi_H v|_K$ is given by
\begin{equation*}
\Pi_H v |_K=\frac{1}{|K|}\int_K v\,dx. 
\end{equation*}
Classical results (see, e.g.,~\cite{Poincareoriginal, Poincare}) state the following (local) stability and approximation property of $\Pi_H$:
\begin{align}
\|\PiH  v\|_{L^2(\elem)}&\leq \|v\|_{L^2(\elem)} &&\text{for all }v\in L^2(\elem),\label{eq:L2stab}\\
\|v-\PiH v\|_{L^2(\elem)}&\leq {\pi}^{-1}H \|\nabla v\|_{L^2(\elem)} &&\text{for all }v\in H^1(\elem).\label{eq:L2approx}
\end{align}

\subsection{Prototypical approximation}
\label{sec:protapp}

Similarly as in~\cite{Owh17,Peterseim2021}, we henceforth derive prototypical problem-adapted ansatz spaces with uniform approximation rates independent of the oscillations of the coefficients. The choice of ansatz space can be motivated by a reduced basis approach in the right-hand side. In order to make this clear, we for now consider the parameter-dependent right-hand side $f_\vartheta \in H^s(\Omega)$, $s \in [0,1]$ with~$\vartheta$ denoting the right-hand side parameter. Besides the $H^s$-regularity, $f_\vartheta$ can be a general parameter-dependent function which might depend on $\vartheta$ in a non-smooth and non-affine way. For non-affine right-hand sides, one typically computes an approximate affine decomposition using appropriate interpolation techniques, e.g., the empirical interpolation method (cf.~\cref{rem:appaffdec}). 

Due to the lack of smoothness, sophisticated interpolation methods as the empirical interpolation method are possibly ineffective. Indeed, for $H^s$-regular functions, the $L^2$-projection onto $\TH$-piecewise constants has optimal approximation properties. Thus, we employ the approximate affine decomposition 
\begin{equation*}
\PiH f_\vartheta = \sum_{K \in \T_{H}} (\PiH f_\vartheta)|_K \; \one_K.
\end{equation*}
We use the reduced basis space
\begin{equation}
\label{eq:idealV}
\V_{H,\mu} \coloneqq \linh\{\mathcal A_\mu^{-1}\one_K\with K \in \TH\}
\end{equation}
which is not adapted to the precise parameter dependence of $f_\vartheta$ and therefore works for any $H^s$-regular right-hand side. Note that we can only expect an algebraic decay of the approximation error in the number of basis functions. More precisely, defining the approximation $u_{H,\mu,\vartheta}$ to the solution $u_{\mu,\vartheta} = \mathcal A^{-1}_{\mu}f_\vartheta\in\V$ in the finite dimensional space $\V_{H,\mu}$ as
\begin{equation}
\label{eq:collocation}
u_{H,\mu,\vartheta} \coloneqq \mathcal A_{\mu}^{-1} \PiH f_\vartheta,
\end{equation}
we obtain the following error estimate: 
there exists a constant $C>0$ independent of $H, \mu$, and $\vartheta$ such that, for any $f_\vartheta \in H^s(\Omega)$, $s \in [0,1]$, the following estimate holds
\begin{equation}
\label{eq:erroruH}
\pi \alpha\vnormof{u_{\mu,\vartheta} - u_{H,\mu,\vartheta}} \leq H \tnormf{f_\vartheta-\PiH f_\vartheta}{\Omega} \leq C H^{1+s} |f_\vartheta|_{H^s(\Omega)},
\end{equation}
where $|f_\vartheta|_{H^s(\Omega)}$ denotes the $H^s$-seminorm of $f_\vartheta$. This result can be proved following the same steps as in the proof of~\cite[Lem.~3.1]{HaPe21b}.

Due to the algebraic convergence of the prototypical approximation in~\cref{eq:erroruH}, the dimension of the reduced basis space $\V_{H,\mu}$ is typically large. Additionally, the canonical basis functions~$\{\mathcal{A}_\mu^{-1}\one_K\with K \in \TH\}$ of~\cref{eq:idealV} are non-local and decay very slowly. Hence, without modification, such approaches are intractable in practice. In order to cure this problem, we subsequently introduce a localization approach identifying (almost) local basis functions of~\cref{eq:idealV}. This allows one to reduce the required storage and the computational costs to a practically feasible level.
For the ease of presentation, we henceforth again consider a parameter-independent right-hand side and drop the dependence of  $u_{H,\mu,\vartheta}$ and $u_{\mu,\vartheta}$ on $\vartheta$.
\subsection{Super-localization technique}
As localization approach, we use a variant of the SLOD from~\cite{HaPe21b}. The SLOD is an improvement of the LOD~\cite{MaP14,HeP13} that exhibits super-exponential localization properties.
For the definition of the method, we first introduce some notations.
Given a union of elements $S\subset \Omega$,
the $\ell$-th order patch of $S$ is defined recursively by 
\begin{align*}
\patch^0(S) \coloneqq S,\qquad \patch^{\ell + 1}(S) \coloneqq \bigcup\left\lbrace K\in\TH \with K \cap \patch^\ell(S) \neq \emptyset \right\rbrace.
\end{align*}
In order to simplify the notation in the subsequent derivation, we fix an element $K\in\T_H$ and the oversampling parameter $\ell\in\mathbb{N}$. 
We will refer to the $\ell$-th order patch of $K$ by $\omega \coloneqq \mathsf{N}^\ell(K)$  and make the meaningful assumption that the patch $\omega$ does not coincide with the whole domain $\Omega$. We define $\Sigma \coloneqq \partial \omega \backslash\partial \Omega$ and denote by $\T_{H,\omega}$ the submesh of $\TH$ with elements in $\omega$. 

For the patch $\omega$, the SLOD aims at identifying an $L^2$-normalized source term  $g = g_{K,\ell,\mu} \in \Pnull(\T_{H,\omega})$ that yields a rapidly decaying (or even local) response $\varphi = \varphi_{K,\ell,\mu} \in \V_{H,\mu}$ under the solution operator~$\mathcal A^{-1}_\mu$, i.e., 
\begin{equation*}
\varphi = \mathcal A^{-1}_\mu g.
\end{equation*}
Note that, throughout this paper, we will not distinguish between $\V_\omega$-functions and their  $\V$-conforming extension by zero to the full domain $\Omega$. 
We define a patch-local approximation $\psi = \psi_{K,\ell,\mu}\in \V_\omega$ of the response $\varphi$ by
\begin{equation}
\label{eq:basisfun}
\psi \coloneqq \mathcal A_{\omega,\mu}^{-1}g,
\end{equation} 
where we recall that $\mathcal A_{\omega,\mu}$ denotes the local solution operator corresponding to~\cref{eq:local_pde}.
For a generic right-hand side $g\in  \Pnull(\T_{H,\omega})$, the local basis function $\psi$ is a poor approximation of the ideal basis function $\varphi$. Nevertheless, the appropriate choice of $g$ leads to a small localization error.

\subsection{Choice of optimal local right-hand sides}

Let us recall the definitions and some properties of traces and extensions of $\V_\omega$-functions.  Defining the space $U \coloneqq \V|_\omega$, we denote the trace operator defined for functions in $U$ and restricted to $\Sigma$ by  
\begin{equation}
\label{eq:X}
\mathrm{tr}_\Sigma \colon  U \to X \coloneqq \operatorname{image}\mathrm{tr}_\Sigma \subset  H^{1/2}(\Sigma).
\end{equation}
The space $X$ is a Hilbert space which can be equipped with the norm
\begin{equation}
\label{eq:defnormtracesp}
\|w\|_X \coloneqq \inf\{\vnormf{v}{\omega}\with v \in H^1(\omega),\; \mathrm{tr}_\Sigma v = w\},
\end{equation}
where we recall that $\vnormf{\bullet}{\omega}$ denotes the restriction of the $\V$-norm to $\omega$.
By definition of the $\|\bullet\|_X$-norm, the continuity of the trace operator holds independent of the patch geometry, namely, for all $v \in U$,
\begin{equation}
\label{eq:conttr}
\|\mathrm{tr}_\Sigma v\|_X \leq \vnormf{v}{\omega}.
\end{equation}
Generalizing the textbook result from~\cite[Ch.~III.1,~Eq.~(1.5)]{BrF91}, we can explicitly compute the $X$-norm of a function in $X$ as follows. 
\begin{lemma}[Computation of $\|\cdot\|_X$]
	For any $w \in X$, the weak solution $u_w \in H^1(\omega)$ to the boundary value problem
	\begin{equation}
	\label{eq:pri}
	\left\{
	\begin{aligned}
	-\Delta u_w + u_w &= 0\, && \text{in}\; \omega,\\
	u_w &= w\; &&\text{on}\; \Sigma,\\
	\nabla u_w \cdot \nu  &= 0\; &&\text{on}\; \partial \omega\cap (\Gamma_1\cup \Gamma_2),\\
	u_w &= 0\; &&\text{on}\; \partial \omega \cap \Gamma_3
	\end{aligned}
	\right.
	\end{equation}
	satisfies
	\begin{equation}
	\label{eq:normeq}
	\vnormf{u_w}{\omega} = \|w\|_X.
	\end{equation}
\end{lemma}
\begin{proof}
	For any $w \in X$, the functional 
	$$J\colon \{v \in H^1(\omega)\with \mathrm{tr}_\Sigma v = w\}\to \mathbb R,\quad  v \mapsto \frac{1}{2}\vnormf{v}{\omega}^2$$ 
	is strictly convex. Thus, the condition that the functional's Gateaux derivative is zero for any direction $\eta \in \V_\omega$ is a sufficient optimality condition. For all $\eta \in \V_\omega$, we obtain
	\begin{equation*}
	\Big(\frac{\mathrm{d}}{\mathrm{d}t}J(u_w+t\eta)\Big)\Big\vert_{t = 0} = \vspf{u_w}{\eta}{\omega}= 0,
	\end{equation*}
	which is the weak formulation of~\cref{eq:pri}.
\end{proof}

This result can be used to define a right-inverse of the trace operator, denoted by $\mathrm{tr}_\Sigma^{-1} \colon X \to U$, by setting $\mathrm{tr}_\Sigma^{-1}w \coloneqq u_w$. By definition, the right-inverse is continuous, more specifically, by~\cref{eq:normeq}, there holds, for all $w \in X$,
\begin{equation}
\label{eq:conttrinv}
\vnormf{\mathrm{tr}_\Sigma^{-1}w}{\omega} = \|w\|_X.
\end{equation}
The conormal derivative of $\psi$ at the boundary segment $\Sigma = \partial \Omega \backslash\partial \omega$ is a functional in~$X^\prime$ which is defined for all $w \in X$ as follows
\begin{equation}
\langle a_\mu \nabla \psi \cdot \nu,w\rangle_{X^\prime\times X} = A_{\omega,_\mu}(\psi,\mathrm{tr}^{-1}_\Sigma w) - \tspf{g}{\mathrm{tr}^{-1}_\Sigma w}{\omega},
\end{equation}
where $\langle \cdot,\cdot\rangle_{X^\prime\times X}$ denotes the duality pairing.
Note that, the conormal derivative is independent of the choice of extension operator $\mathrm{tr}_\Sigma^{-1}$.

The following crucial observation establishes a connection between the localization error and the norm of the conormal derivative of $\psi$. For any $v \in \V$, it holds 
\begin{equation}
\label{e:nd1}
A_\mu (\psi-\varphi,v) = A_{\omega,_\mu}(\psi  ,v)-(g,v)_{L^2(\omega)} =  \langle a_\mu \nabla \psi \cdot \nu,\mathrm{tr}_\Sigma\, v|_\omega\rangle_{X^\prime\times X}.
\end{equation}
It is possible to explicitly compute the $X^\prime$-norm of a functional by generalizing the textbook result~\cite[Ch.~III.1,~Eq.~(1.8)]{BrF91} as follows.

\begin{lemma}[Computation of $\|\cdot\|_{X^\prime}$]
	\label{lemma:xprimenorm}
	For any $q \in X^\prime$, the weak solution $u_q \in H^1(\omega)$ to the boundary value problem
	\begin{equation}
	\label{eq:prxp}
	\left\{
	\begin{aligned}
	-\Delta u_q + u_q &= 0\, && \text{in}\; \omega,\\
	\nabla u_q \cdot \nu &= q\; &&\text{on}\; \Sigma,\\
	\nabla u_q \cdot \nu  &= 0\; &&\text{on}\; \partial \omega\cap (\Gamma_1\cup \Gamma_2),\\
	u_q &= 0\; &&\text{on}\; \partial \omega \cap \Gamma_3
	\end{aligned}
	\right.
	\end{equation}
	satisfies
	\begin{equation}
	\label{eq:normeqXprime}
	\vnormf{u_q}{\omega} = \|q\|_{X^\prime}.
	\end{equation}
\end{lemma}
\begin{proof}
	Employing that $u_q \in U$ is the weak solution to~\cref{eq:prxp} and using~\cref{eq:defnormtracesp,eq:normeq}, we obtain
	\begin{align*}
	\vnormf{u_q}{\omega} &= \sup_{v \in U} \frac{\vspf{u_q}{v}{\omega}}{\vnormf{v}{\omega}} = \sup_{v \in U} \frac{\langle q, \mathrm{tr}_\Sigma v\rangle_{X^\prime \times X} }{\vnormf{v}{\omega}}=\sup_{w \in X} \frac{\langle q, \mathrm{tr}_\Sigma u_w\rangle_{X^\prime \times X} }{\vnormf{u_w}{\omega}}\\& = \sup_{w \in X} \frac{\langle q, w\rangle_{X^\prime \times X} }{\|w\|_{X}} = \|q\|_{X^\prime}
	\end{align*}
	which proves the desired result.
\end{proof}
For the sake of notation, we abbreviate the mapping of a right-hand side $g$ to the element $u_q$ that solves~\cref{eq:prxp} for $q = a_\mu  \nabla \psi\cdot \nu \in X^\prime$ with $\psi$ being the basis function corresponding to $g$, by the operator $\mathcal R = \mathcal R_{K,\ell,\mu}$, i.e., 
\begin{equation}
\label{eq:R}
\mathcal R\colon \Pnull(\T_{H,\omega})\to U, \; g \mapsto u_q.
\end{equation}
By observation~\cref{e:nd1} and~\cref{lemma:xprimenorm}, we achieve a minimal localization error by choosing the right-hand side as the solution to the energy minimization problem
\begin{equation}
\label{eq:chocieg}
g = \argmin_{p \in \Pnull(\T_{H,\omega})}\frac{\vnormf{\mathcal R p}{\omega}^2}{\tnormf{p}{\omega}^2}.
\end{equation}
It shall be noted that this problem is equivalent to solving an eigenvalue problem as the following remark states. 
\begin{remark}[Equivalent eigenvalue problem]
	\label{rem:eqevp}
	Instead of solving~\cref{eq:chocieg}, it is equivalent to compute the eigenvector corresponding to the smallest eigenvalue of the following generalized eigenvalue problem 
	\begin{equation}
	\label{eq:eqevp}
	Cx = \lambda Dx
	\end{equation}
	with matrices $C,D \in \mathbb R^{J\times J}$, $J \coloneqq \# \T_{H,\omega}$, defined as
	\begin{align*}
	C_{ij} = \vspf{\R \one_{T_j}}{\R \one_{T_i}}{\omega},\; D_{ij} \coloneqq \tspf{\one_{T_j}}{\one_{T_i}}{\omega},
	\end{align*}
	with $\{T_j\with j = 1,\dots J\}$ being some numbering of the elements in $\T_{H,\omega}$.
\end{remark}


\subsection{Localization error indicator}

We introduce the local error indicator $\sigma_{K,\mu} = \linebreak \sigma_{K,\mu}(H,\ell) \coloneqq \vnormf{\mathcal R g}{\omega}$ which coincides with the $X^\prime$-norm of the conormal derivative of~$\psi$. Taking the maximum, we obtain an error indicator for the method's localization error
\begin{equation}
\label{eq:sigma}
\sigma_\mu = \sigma_\mu(H,\ell) \coloneqq \max_{K\in\TH} \sigma_{K,\mu}.
\end{equation}
Based on extensive numerical studies in~\cite{HaPe21b,Freese-Hauck-Peterseim,BFP22} and a justification~\cite[Thm.~7.3]{HaPe21b} that relies on a conjecture from spectral geometry on the decay of Steklov eigenfunctions, we conjecture that $\sigma$ decays super-exponentially in the oversampling parameter~$\ell$.

\begin{conjecture}[Super-exponential decay of $\sigma_\mu$]\label{conj:sexpdec}
	For all parameters $\mu \in \pset$, the quantity~$\sigma_\mu$ decays super-exponentially in $\ell$, i.e., there exists a constants $C_\sigma(H,\ell)>0$ depending  polynomially on $H$ and $\ell$ and independent of $\mu$ and $C>0$ independent of $H,\ell,$ and $\mu$ such that
	\begin{equation*}
	\sigma_\mu(H,\ell) \leq C_\sigma(H,\ell)\exp(-C\ell^{\frac{d}{d-1}}).
	\end{equation*}
\end{conjecture}
Note that, using LOD techniques, one can prove a pessimistic bound still guaranteeing an exponential decay of $\sigma$ as $\ell$ is increased, cf.~\cite[Lem.~6.4]{HaPe21b}.

\subsection{Super-localized multiscale method}

For computing a discrete approximation of $u_\mu \in \V$ in the space
\begin{equation*}
\V_{H,\ell,\mu} \coloneqq \linh\{\psi_{K,\ell,\,u}\with K \in \TH\},
\end{equation*}
we employ the so-called collocation variant of the SLOD~\cite[Rem.~5.1]{HaPe21b}. This variant has the advantage that, compared to computing the Galerkin solution in $\V_{H,\ell,\mu}$, no inner products between basis functions need to be computed. Instead, we define the discrete approximation $u_{H,\ell,\mu}\in \V_{H,\ell,\mu}$ as
\begin{equation}
\label{eq:colsol}
u_{H,\ell,\mu} = \sum_{K \in \mathcal T_H}c_K \, \psi_{K,\ell,\mu},
\end{equation}
where $(c_K)_{K \in \TH}$ are the coefficients of the expansion of $\PiH f$ in terms of the basis functions $g_{K,\ell,\mu}$
\begin{equation*}
\PiH f =\sum_{K\in\TH} c_K g_{K,\ell,\mu}.
\end{equation*}

\begin{remark}[Connection to isogeometric analysis]
	In 1d, the SLOD basis function for the diffusion problem with constant coefficient $a_\mu$ coincides with the quadratic B-splines. We refer to~\cite[Fig.~4.1]{HaPe21b} for an illustration. 
\end{remark}

\section{Reduced basis super-localized orthogonal decomposition}
\label{sec:RB-SLOD}

This section introduces the reduced basis super-localized orthogonal decomposition (RB-SLOD) method which combines the SLOD presented in the last section with a reduced basis approach enabling a significant acceleration in the context of parametric problems. The algorithm consists of an offline phase and an online phase. The offline phase is executed only once and performs precomputations which are employed in the online phase for the rapid computation of an approximation to the solution $u_\mu$.

\subsection{Offline phase}

Henceforth, we fix an element $K\in \TH$ and its $\ell$-th order patch $\omega \coloneqq \mathsf N^\ell(K)$. Given a parameter $\mu\in\pset$, the main computational effort required for the computation of the basis functions $\psi\in \V_{H,\mu}$ in~\cref{eq:basisfun} and $g\in \Pnull(\TH)$ in~\cref{eq:chocieg} is the computation of the set
\begin{equation*}
\{ \chi_{T,\mu}\with T \in \T_{H,\omega}\}
\end{equation*}
with 
\begin{equation}
\label{eq:chi}
\chi_{T,\mu} \coloneqq \mathcal A_{\omega,\mu}^{-1}\one_T.
\end{equation}
For all $T \in \T_{H,\omega}$, we employ a reduced basis approach for deriving approximations of $\chi_{T,\mu}$ that can be rapidly evaluated. Henceforth, we also fix the element $T \subset \mathsf N^\ell(K)$ and omit the subscript $T$ for the ease of readability.

\subsubsection{Initialization}

The first step in the reduced basis approach is the selection of a training set of parameters $\ptr \subset \pset$ with size given by the prescribed integer $L$. Possible options to select the elements of $\ptr$ are random sampling or structured grids.~\cref{algo:initialization1} shows the initialization of the set of reduced basis parameters~$\prb$ and the set of reduced basis functions $\mathcal W^\mathrm{rb}$.

\begin{algorithm}[h]
	\begin{algorithmic}[1]
		\State $\ptr \gets  \{\mu_1,\dots,\mu_L\}$ with selected parameters from $\pset$
		\State $\prb \gets \{\mu_1\}$
		\State $\mathcal W^\mathrm{rb} \gets \{\chi_{\mu_1}\}$ with $\chi_{\mu_1} \coloneqq \mathcal{A}_{\omega,\mu_1}^{-1} \one_T$
	\end{algorithmic}
	\caption{(Offline -- Initialization with starting parameters).}
	\label{algo:initialization1}
\end{algorithm}

\subsubsection{Error estimation}

Assume that we have already selected a set of reduced basis parameters $\prb = \{\mu_1,\dots,\mu_M\}$ with corresponding reduced basis functions $\mathcal W^\mathrm{rb} = \{\chi_{\mu_1},\dots,\chi_{\mu_M}\}$. Given the new parameter $\mu \in \ptr \backslash \prb$, we aim at
\begin{enumerate}
	\item computing an approximation $\chi^\mathrm{rb}_\mu$ of $\chi_{\mu} = \mathcal A^{-1}_{\omega,\mu}\one_T$ in the span of $\mathcal W^\mathrm{rb}$ and
	\item efficiently estimating the approximation error.
\end{enumerate}
To perform the first task, we compute the Galerkin projection onto the span of $\mathcal W^\mathrm{rb}$, i.e., we seek $\chi_{\mu}^\mathrm{rb}=\sum_{m=1}^M c_m \chi_{\mu_m} \in \operatorname{span}\mathcal W^\mathrm{rb}$ such that, for all $w \in \operatorname{span}\mathcal W^\mathrm{rb}$,
\begin{equation}
\label{eq:rbbestapprox}
A_{\omega,\mu}(\chi_{\mu}^\mathrm{rb},w) = A_{\omega,\mu}(\chi_{\mu},w) = \tspf{\one_T}{w}{\omega}.
\end{equation}
Choosing test functions in $\mathcal W^\mathrm{rb}$,~\cref{eq:rbbestapprox} turns into a small linear system with $(c_m)_{m=1}^M$ being the unknowns. Note that $\chi_{\mu}$ itself is not required for computing the approximation. 

For estimating the approximation error, we compute the Riesz representation $\tau_\mu$ of the residual functional, i.e., we seek $\tau_\mu \in \V_\omega$ such that, for all $w \in \V_\omega$,
\begin{equation}
\label{eq:defr}
\vspf{\tau_\mu}{w}{\omega} =  A_{\omega,\mu}(\chi_{\mu}^\mathrm{rb},w) - \tspf{\one_T}{w}{\omega}
\end{equation}
and define the following a posteriori error estimator
\begin{equation}
\label{eq:Delta}
\Delta_\mu \coloneqq \vnormf{\tau_\mu}{\omega}.	
\end{equation}
The following lemma proves that this error estimator is reliable and efficient. 
\begin{lemma}[Reliability and efficiency of estimator]
	\label{lemma:estimator} 
	For any given $\mu \in \ptr \backslash \prb$, the reduced basis error for $\chi_\mu$ can be bounded from below and above by $\Delta_\mu$ defined in~\cref{eq:defr}, i.e., there holds
	\begin{equation*}
	\beta^{-1}\Delta_\mu\leq \vnormf{\chi_{\mu}-\chi_{\mu}^\mathrm{rb}}{\omega} \leq \alpha^{-1}\Delta_\mu.
	\end{equation*}
\end{lemma}
\begin{proof}
	We define $e \coloneqq \chi_{\mu}-\chi_{\mu}^\mathrm{rb}$. Using~\cref{eq:conditions4laxmilgram,eq:defr} together with $\V_\omega \subset \V$, we obtain
	\begin{equation*}
	\vnormf{\tau_\mu}{\omega}^2 = A_{\omega,\mu}(\chi_{\mu}^\mathrm{rb},\tau_\mu) - \tspf{\one_T}{\tau_\mu}{\omega}=A_{\omega,\mu}(e,\tau_\mu)\leq \beta \vnormf{e}{\omega} \vnormf{\tau_\mu}{\omega}
	\end{equation*}
	which, after dividing by $\vnormf{\tau_\mu}{\omega}$, implies the first inequality. For the second inequality, we similarly obtain
	\begin{equation*}
	\alpha \vnormf{e}{\omega}^2 \leq A_{\omega,\mu}(e,e) = \tspf{\one_T}{e}{\omega}- A_{\omega,\mu}(\chi_{\mu}^\mathrm{rb},e) \leq \vnormf{\tau_\mu}{\omega}\vnormf{e}{\omega}.\qedhere
	\end{equation*}
\end{proof}

Using the affine decomposition~\cref{eq:decA}, and recalling that $\chi_{\mu}^\mathrm{rb}=\sum_{m=1}^M c_m \chi_{\mu_m}$, the a posteriori error estimator can be efficiently computed observing that
\begin{equation}
\label{eq:tau}
\tau_\mu = \sum_{q=1}^Q \theta_q(\mu)\sum_{m = 1}^M c_m s_{q,m} - p,
\end{equation}
where the $s_{q,m} \in \V_\omega$ are the unique functions satisfying, for all $w \in \V_\omega$,
\begin{equation}
\label{eq:sqk}
\vspf{s_{q,m}}{w}{\omega} = B_{\omega,q}(\chi_{T,\mu_m},w)
\end{equation}
and $p \in \V_\omega$ is the unique function satisfying, for all $w \in \V_\omega$,
\begin{equation}
\label{eq:p}
\vspf{p}{w}{\omega} = \tspf{\one_T}{w}{\omega}.
\end{equation}

Denoting, the sets in which we store the functions $\{s_{q,m}\with m = 1,\dots,M\}$ by $\mathcal S_q$, we can initialize the Riesz representatives as shown in~\cref{algo:initialization2}.
\begin{algorithm}[h]
	\begin{algorithmic}[1]
		\For{$q = 1,\dots,Q$}
		\State $\mathcal S_q \gets \{s_{q,1}\}$ with $s_{q,1}$ solving~\cref{eq:sqk}
		\EndFor
		\State compute $p$ by~\cref{eq:p}
	\end{algorithmic}
	\caption{(Offline -- Initialization Riesz representatives).}
	\label{algo:initialization2}
\end{algorithm}

\subsubsection{Greedy search}
The greedy search algorithm iterates through the parameters $\mu \in  \ptr\backslash\prb$ and, for each $\mu$, estimates the error when approximating $\chi_{\mu}$ by $\chi_{\mu}^\mathrm{rb}$ using the error estimator~$\Delta_{\mu}$ from~\cref{eq:Delta}. It then selects the parameter $\mu$ for which the estimator is largest and adds it to the set $\prb$ . The corresponding function $\chi_{\mu}$ is then computed by~\cref{eq:chi} and added to the set of reduced basis functions $\mathcal W^\mathrm{rb}$. For numerical stability reasons, we additionally perform an orthogonalization with respect to $\vspf{\cdot}{\cdot}{\omega}$ using a Gram--Schmidt-type algorithm. Note that, due to the typically small size of $\prb$ the numerical instability issues of the Gram--Schmidt algorithm are not noticeable in practice. 

Given a tolerance $\mathsf{tol}>0$, this procedure is repeated until the training error satisfies the following stopping criterion
\begin{equation}
\label{eq:esterr}
\frac{\max_{\mu \in \ptr} \Delta_\mu }{ \vnormf{p}{\omega}}
\leq \mathsf{tol},
\end{equation}
with $p$ defined in~\cref{eq:p}. It is justified to use $p$ instead of $\chi_{\mu}$ in the denominator of~\cref{eq:esterr} as
\begin{equation*}
\vnormf{p}{\omega} = \|\one_T\|_{\V^\prime_\omega} \approx \vnormf{ \chi_{\mu}}{\omega},
\end{equation*}
where $\|\cdot\|_{\V^\prime_\omega}$ denotes the canonical norm on the dual space $\V^\prime_\omega$. In particular, this means that  both quantities have the same scaling in $H$.~\cref{algo:greedysearch} shows an implementation of the greedy search in pseudo code.

\begin{algorithm}[h]
	\begin{algorithmic}[1]
		\While{$\frac{\underset{\mu \in \ptr}{\operatorname{max}} \Delta_{\mu}}{\vnormf{ p}{\omega}}>\mathsf{tol}$}
		\State $M \gets \#\prb$
		\State // compute error estimators
		\For{$\mu \in \ptr\backslash \prb$}
		\State compute $\chi_{\mu}^\mathrm{rb}$ by~\cref{eq:rbbestapprox} \label{line:chirb}
		\State compute $\tau_\mu$ by~\cref{eq:tau}
		\State $\Delta_\mu \gets \vnormf{ \tau_\mu}{\omega}$ by~\cref{eq:Delta}
		\EndFor
		\State // select element with largest estimator
		\State $\mu_{M+1} \gets \underset{ \mu \in \ptr \backslash \prb }{\operatorname{argmax}} \Delta_\mu$
		\State $\prb \gets \prb \cup \mu_{M+1}$
		\State compute $\chi_{\mu_{M+1}}$ by~\cref{eq:chi}\label{line:chimu}
		\State $\mathcal W^\mathrm{rb}\gets \mathcal W^\mathrm{rb} \cup \{\chi_{\mu_{M+1}}\}$
		\State // update Riesz representatives
		\For{$q = 1,\dots,Q$}
		\State $\mathcal S_q \gets \mathcal S_q \cup \{s_{q,M+1}\}$ with $s_{q,M+1}$ by~\cref{eq:sqk}
		\EndFor
		\EndWhile
	\end{algorithmic}
	\caption{(Offline -- Greedy search).}
	\label{algo:greedysearch}
\end{algorithm}

The computation of $\chi_{\mu}^\mathrm{rb}$ in Line~\ref{line:chirb} of~\cref{algo:greedysearch} can be accelerated employing the affine decomposition~\cref{eq:decA} of the bilinear form. 
After selecting an element~$\mu_{M+1}$, adding it to~$\prb$ and computing the function $\chi_{\mu_{M+1}}$, one may precompute the inner products 
\begin{equation*}
\tspf{\one_T}{\chi_{\mu_{M+1}}}{\omega}\quad \text{and}\quad  B_{\omega,q}(\chi_{\mu_{m}},\chi_{\mu_{M+1}}),  \quad q = 1,\dots,Q,\;m = 1,\dots,M.
\end{equation*}
This is especially important, as the training set $\ptr$ is typically large and the computation of~\cref{eq:rbbestapprox} in Line~\ref{line:chirb} has to be repeated many times. By precomputation, one can reduce the complexity of Line~\ref{line:chirb} to polynomial complexity in $\#\T_{H,\omega}$. 

\begin{remark}[Decay of training error]
	\label{rem:decdelta}
	In practice, provided that the bilinear form $A_\mu$ depends smoothly on the parameter~$\mu$, typically a (sub-)exponential decay of the training error~\cref{eq:esterr} can be observed. For a numerical investigation, see the experiments in~\cref{sec:numexp}. 
	Theoretically, it is difficult to derive explicit rates of decay of the training error, see e.g.~\cite{BMP12,AbH15,ohlberger2016reduced} for an investigation. 
\end{remark}

All the computations described in this sections need to be repeated for all $T \in \T_{H,\omega}$ and $K\in\TH$. This can be done in parallel. Hence, the final output of the offline phase is the collection of RB spaces $\mathcal W^{\textrm{rb}}$ as the element $K$ varies in $\TH$ and $T$ varies in $\T_{H,\omega}$. 

\subsection{Online phase}

Given any parameter value of interest $\mu \in \pset$, the online phase rapidly computes approximate SLOD basis functions using the precomputations that were performed in the offline phase. The RB-SLOD solution is then obtained after solving a sparse coarse linear system. 

\subsubsection{Basis computation}

As before, we fix the element $K \in \TH$ and denote its $\ell$-th order patch by $\omega \coloneqq \mathsf N^\ell(K)$. We aim at constructing an approximation of the SLOD basis function $\psi = \psi_{K,\ell,\mu}$ from~\cref{eq:basisfun} by 
\begin{equation}
\label{eq:psirb}
\psi^\mathrm{rb} \coloneqq \sum_{T \in \T_{H,\omega}} c_T\chi_{T,\mu}^\mathrm{rb},
\end{equation}
where the coefficients $(c_T)_{T \in \T_{H,\omega}}$ will be determined subsequently and $\chi_{T,\mu}^\mathrm{rb}$ denotes the RB approximation of  $ \chi_{T,\mu} = \mathcal{A}_{\omega,\mu}^{-1}\one_T$  computed by solving the small linear system in~\cref{eq:rbbestapprox}. Note that, unlike the SLOD basis functions~\cref{eq:basisfun}, their RB approximations~\cref{eq:psirb} admit no $L^2$-regular right-hand side via the PDE operator $\mathcal A_{\omega,\mu}$. Nevertheless, by analogy with the super-localization technique for non-parametric PDEs, we define $g^\mathrm{rb} = g^\mathrm{rb}_{K,\ell,\mu}$ as
\begin{equation}
\label{eq:grb}
g^\mathrm{rb} \coloneqq \sum_{T \in \T_{H,\omega}} c_T \one_T,
\end{equation}
where the coefficients $(c_T)_{T \in \T_{H,\omega}}$ are the same as in~\cref{eq:psirb}. In particular, we have established an one-to-one relationship between $g^\mathrm{rb}\in \Pnull(\T_{H,\omega})$ and $\psi^\mathrm{rb}\in H^1(\omega)$. Unfortunately, the straightforward adaptation of the super-localization procedure illustrated in~\cref{sec:SLOD} can not be pursued here, since the conormal derivative of $\psi^\mathrm{rb}$, in general, does not exist as a functional in $X^\prime$ with $X$  defined in~\cref{eq:X}. Recalling the definition of the extension operator $\mathrm{tr}^{-1}_\Sigma\colon  U \to X$ with $\Sigma = \partial \omega \backslash \partial \Omega$ and $U = \V|_\omega$ in~\cref{eq:X}, we define an analogue to the conormal derivative of $\psi^\mathrm{rb}$ similarly as~\cref{eq:conormalder} by
\begin{equation}
\label{eq:conormalder}
\langle a_{\mu}\nabla \psi^\mathrm{rb} \cdot \nu, w\rangle_{X^\prime\times X} \coloneqq A_{\omega,\mu}(\psi^\mathrm{rb} ,\mathrm{tr}_\Sigma^{-1} w) - (g^\mathrm{rb},\mathrm{tr}_\Sigma^{-1}w)_{L^2(\omega)}.
\end{equation}
Note that this functional clearly is an element of $X^\prime$ which however depends on the choice of extension operator, i.e., for two different extension operators, the respective functionals might not coincide. Nevertheless, the $X^\prime$-norm of the difference of the functionals is determined by the tolerance~\cref{eq:esterr} of the greedy search algorithm, i.e., for small tolerances the functionals corresponding to different extension operators (almost) coincide.

We denote with $\mathcal R^\mathrm{rb} = \mathcal R_{K,\ell,\mu}^\mathrm{rb}\colon \Pnull(\T_{H,\omega})\to H^1(\omega)$ the operator mapping a right-hand side $g^\mathrm{rb}$ as in~\cref{eq:grb} to $v\in H^1(\omega)$ being the unique solution to~\cref{eq:prxp} with $q = a_\mu \nabla \psi^\mathrm{rb}\cdot \nu\in X^\prime$.
Similar as in~\cref{eq:chocieg}, we choose the right-hand side $g^\mathrm{rb}$ so that the conormal derivative of the basis function $\psi^\mathrm{rb}$ is minimal which translates to a small localization error. More precisely, we choose
\begin{equation}
\label{eq:chociegrb}
g^\mathrm{rb} \coloneqq \argmin_{p \in \Pnull(\T_{H,\omega})}\frac{ \vnormf{\mathcal R^\mathrm{rb}p}{\omega}^2}{\tnormf{p}{\omega}^2}.
\end{equation}


In~\cref{algo:basiscomp}, we denote by $\B$ and $\mathcal G$ the sets, where the RB-SLOD basis functions $\{\psi^\mathrm{rb}_{K,\ell,\mu},\, K\in\TH,\, \mu\in\ptr\}$ and the corresponding $\{g^\mathrm{rb}_{K,\ell,\mu},\, K\in\TH,\, \mu\in\ptr\}$ are stored, respectively. See \cref{fig:slodbasis} for a depiction of $g^\mathrm{rb}$ and $\psi^\mathrm{rb}$ for a parametric reaction-convection-diffusion problem for several choices of parameters.

\begin{algorithm}[h]
	\begin{algorithmic}[1]
		\State $\mathcal B,\, \mathcal G \gets  \{\}$
		\For{$K \in \TH$}
		\State // compute reduced basis approximations
		\For{$T \in \T_{H,\,\mathsf N^\ell(K)}$}
		\State compute $\chi_{T,\mu}^\mathrm{rb}$ by~\cref{eq:rbbestapprox} \label{line:computerbchar}
		\EndFor	
		\State // compute and save basis functions
		\State compute $g_{K,\ell,\mu}^\mathrm{rb}$ by~\cref{eq:chociegrb}\label{line:g}
		\State $\mathcal G \gets \mathcal G \cup \{g_{K,\ell,\mu}^\mathrm{rb}\}$
		\State obtain $\psi_{K,\ell,\mu}^\mathrm{rb}$ by~\cref{eq:psirb} with coefficients from~\cref{eq:grb}\label{line:computeslodbasis}
		\State $\mathcal B \gets \mathcal B \cup \{\psi_{K,\ell,\mu}^\mathrm{rb}\}$
		\EndFor		
	\end{algorithmic}
	\caption{(Online -- RB-SLOD Basis computation).}
	\label{algo:basiscomp}
\end{algorithm}

The above computations need to be repeated for all $K\in \TH$. The intermediate output of the online phase is the pair of sets $\mathcal B$ and $\mathcal G$, collecting all the basis functions and corresponding right-hand sides. This information will be then used during the \emph{coarse solve} for the computation of the RB-SLOD approximation.

\subsubsection{Coarse solve}

For the computation of the RB-SLOD approximation of $u_\mu$, we proceed similarly as in~\cref{eq:colsol}. First, we express $\PiH f$ as linear combination of the right-hand sides in $\mathcal G$, namely we compute coefficients $(c_K)_{K\in\TH}$ such that
\begin{equation}
\label{eq:coeffs}
\PiH f = \sum_{K \in \T_{H}} c_K g^\mathrm{rb}_{K,\ell,\mu},
\end{equation}
cf.~\Cref{algo:coarsesolve}. Then, following the collocation approach, we define the RB-SLOD solution as the linear combination of the basis functions in $\B$ with exactly these coefficients, i.e., 
\begin{equation}
\label{eq:urbslod}
u_{H,\ell,\mu}^\mathrm{rb} \coloneqq \sum_{K \in \T_{H}}c_K \psi_{K,\ell,\mu}^\mathrm{rb}.
\end{equation}

In~\cref{algo:coarsesolve}, we denote by $\mathbf G \in \mathbb{R}^{\#\TH\times \#\TH}$ the matrix which has the ($\TH$-piecewise constants) functions in $\mathcal G$ as its columns, and by $\mathbf f$ the load vector having the element values of $\PiH f$ as its entries. Note that in~\cref{algo:coarsesolve}, no inner product between basis functions are computed.

\begin{algorithm}[h]
	\begin{algorithmic}[1]
		\State assemble $\mathbf G$ from $\mathcal G$
		\State compute $\mathbf f$
		\State solve $\mathbf c = \mathbf G\backslash \mathbf f$
		\State compute $u_{H,\ell,\mu}^\mathrm{rb}$ by~\cref{eq:urbslod} with coefficients $(\mathbf c_i)_{i=1,\dots,\#\TH}$\label{line:addupbasisfun}
	\end{algorithmic}
	\caption{(Online -- Coarse solve)}
	\label{algo:coarsesolve}
\end{algorithm}

\begin{remark}[Built-in reduced basis approach in the right-hand side]
	\label{rem:rbinrhs}
	It is noteworthy that the right-hand side $f$ first appears in~\cref{algo:coarsesolve}, i.e., the offline phase and the basis computation are completely independent of $f$. This feature makes the RB-SLOD method suitable for applications, where the solution for various right-hand sides is of interest. In particular, parameter-dependent and possible non-affine right-hand sides do not pose any further difficulties. In the implementation, the right-hand side parameters only concern the coarse solve and the actual reduced basis approach is independent of these parameters, see also~\cref{sec:protapp} and the numerical example in~\cref{sec:numexp2}.
\end{remark}

\subsection{Practical implementation}
In a practical implementation, we perform a fine-scale discretization of the infinite dimensional patch problems~\cref{eq:chi,eq:defr,eq:sqk,eq:p}. For any patch $\omega = \mathsf N^\ell(K)$, we obtain a fine patch mesh $\T_{h,\omega}$ by uniform refinement of the coarse patch mesh $\T_{H,\omega}$. For discretizing the patch problems, one may employ first order continuous finite elements on the fine patch mesh.
\subsubsection{Practical basis computation}\label{sec:implementation}
Similar to~\cref{eq:psirb}, we define the discrete localized basis function $\psi_{h}^\mathrm{rb}$ supported on the patch $\omega$, as sum of discretized reduced basis approximations $\{\chi_{T,\mu,h}^\mathrm{rb}\with T \in \T_{H,\omega}\}$ with coefficients being the element values of the function $g^\mathrm{rb}_h \in \mathbb P^0(\T_{H,\omega})$. Exploiting that the conormal derivative of $\psi_{h}^\mathrm{rb}$ exists in the $L^2(\partial\omega)$-sense enables alternatives choices of the right-hand sides with reduced computational costs. 

For example, instead of~\cref{eq:chociegrb}, one may minimize the $L^2(\Sigma)$-norm of the conormal derivative $(a_\mu\nabla \psi_{h}^\mathrm{rb}\cdot \nu)|_{\Sigma}$. It is also possible to omit the diffusion tensor in~\cref{eq:choicegh}, i.e., we minimize the normal derivative $(\nabla \psi_{h}^\mathrm{rb}\cdot \nu)|_{\Sigma}$ instead of the conormal derivative. This is justified since by~\cref{eq:unifboundsA}, a small norm of the normal derivative implies a small norm of the conormal derivative. Hence, in the practical implementation, we choose $g^\mathrm{rb}_h$ as the solution to the following constraint minimization problem
\begin{equation}
\label{eq:choicegh}
g^\mathrm{rb}_h \coloneqq \argmin_{p \in \Pnull(\T_{H,\omega})} \frac{\tnormf{\mathcal R^\mathrm{rb}_h p}{\Sigma}^2}{\tnormf{p}{\omega}^2},
\end{equation}
where $\mathcal R^\mathrm{rb}_h\colon \Pnull(\T_{H,\omega})\to L^2(\Sigma)$ denotes the mapping of $g_h^\mathrm{rb}$ to the normal derivative $(\nabla \psi_{h}^\mathrm{rb}\cdot \nu)|_{\Sigma}$ of the corresponding basis function. Such an approach is for example employed in~\cite{BFP22} for the convection-dominated diffusion problem. 

In a practical implementation, instead of~\cref{eq:chociegrb}, we solve the equivalent generalized eigenvalue problem~$Cx = \lambda Dx$ stated in \cref{rem:eqevp} with $\R$ replaced by $\R^\mathrm{rb}_h$. An accelerated assembly of the matrix $C$ can be achieved by precomputing the following four dimensional array in the offline phase  
\begin{align*}
\tspf{\nabla \chi_{T_j,\mu_m}\cdot \nu}{\nabla \chi_{T_i,\mu_l}\cdot \nu}{\Sigma},\qquad l=1,\dots,M_i,\; m = 1,\dots,M_j,\;i,j=1,\dots,J
\end{align*}
with $M_i$ and $M_j$ denoting the cardinality of the parameter sets $\prb$ corresponding to elements $T_i$ and $T_j$, respectively and $J = \#\T_{H,\omega}$. As the mass matrix $D$ is parameter-independent, it only needs to be assembled once.

\subsubsection{Choice of basis}	\label{rem:choiceofbasis}
Note that it is neither guaranteed  that the minimization problems~\cref{eq:chocieg,eq:chociegrb,eq:choicegh} have a unique solution, nor that the right-hand sides in $\B$ form a stable basis of $\Pnull(\TH)$ which is necessary for~\cref{eq:coeffs}. 
In~\cite[App.~B]{HaPe21b}, a stable and efficient algorithm for the selection of the right-hand sides is proposed which can also be applied in the current setting. Based on the observation that selection issues only occur for patches close to the boundary, the algorithm conducts a special treatment of such troubled patches for resolving possible uniqueness and stability issues. For details and a detailed illustrative description of the algorithm, see~\cite[App.~B]{HaPe21b}. 

\subsubsection{Complexity}
Due to the precomputations in the offline phase, no (local or global) fine-scale solves are needed in the online phase. More specifically, the only fine-scale operations in the online phase are (i) adding up (local) fine-scale functions in the computation of the SLOD basis (see~\Cref{algo:basiscomp}, Line~\ref{line:computerbchar} and~\ref{line:computeslodbasis}) and (ii) adding up the SLOD basis functions using the coefficients computed by solving a coarse system (see~\Cref{algo:coarsesolve}, Line~\ref{line:addupbasisfun}).

Note that, if one is only interested in coarse quantities of the solution, fine-scale operations can be completely avoided in the online phase. In the offline phase, one can additionally compute the element averages of the fine-scale functions in~\cref{algo:greedysearch}, Line~\ref{line:chimu}. This can then be used in the online phase for the computing the element averages of the SLOD basis functions without fine-scale operations. By adding up the element averages of the SLOD basis functions using the same coefficients as in~\cref{algo:coarsesolve}, Line~\ref{line:addupbasisfun}, one obtains a $\TH$-piecewise constant approximation of the solution without performing any fine-scale operations in the online phase.

\section{Error analysis}
\label{sec:error_analysis}
In this section, we derive upper bounds for the RB-SLOD approximation error.
As observed in~\cref{rem:choiceofbasis}, it is crucial that the right-hand sides in $\mathcal G$ form a stable basis of~$\Pnull(\TH)$. In practice, this condition can be enforced numerically by the algorithm proposed in~\cite[App.~B]{HaPe21b}. Since, we cannot a priori control the stability of the basis constructed by this algorithm, we make the following assumption for the error analysis of the RB-SLOD.

\begin{assumption}[Riesz stability]\label{a:Rieszbasis}
	\label{ass:stabg}
	The set $\{g_{K,\ell,\mu}^\mathrm{rb}\with K \in \TH\}$ is a Riesz basis of $\Pnull(\mathcal T_H)$, i.e., for all $\mu \in \ptr$, there exists $C_{\mu}(H,\ell)>0$ depending polynomially on $H, \ell$ such that, for all $(c_{K})_{K \in \TH}$,
	\begin{equation*}
	C^{-1}_{\mu}(H,\ell)\sum_{K \in \mathcal T_H} c_{K}^2  \leq \Big\lVert \sum_{K \in \mathcal T_H} c_{K} g_{K,\ell,\mu}^\mathrm{rb}\Bigl\rVert_{L^2(\Omega)}^2 \;\leq C_{\mu}(H,\ell) \sum_{K \in \mathcal T_H} c_{K}^2.
	\end{equation*}
\end{assumption}
We obtain the following error estimate.

\begin{theorem}[Upper bound on the RB-SLOD error]
	\label{th:error}
	Given the parameter $\mu \in \ptr$, we assume that $\{g_{K,\ell,\mu}\with K \in \TH\}$ is stable in the sense of~\cref{a:Rieszbasis}. Let $u_\mu\in\V$ and $u_{H,\ell,\mu}^\mathrm{rb} \in \V_{H,\ell,\mu}^\mathrm{rb}$ denote the unique solutions to~\cref{eq:weak,eq:urbslod}, respectively.
	Then, there exist $C,\,C^\prime>0$ independent of $H,\ell,$ and $\mu$ such that, for all $f \in H^s(\Omega)$, $s \in [0,1]$, 
	\begin{align}
	&\vnormof{u_\mu - u_{H,\ell,\mu}^\mathrm{rb}}\notag\\*
	\begin{split}
	&\qquad  \leq C\hphantom{^\prime} \big(H \tnormf{f-\PiH f}{\Omega} + C_{\mu}^{1/2}(H,\ell)\ell^{d/2}(\sigma_\mu(H,\ell) + \mathsf{tol}  \ell^{d/2})\tnormf{f}{\Omega}\big)\\*
	& \qquad  \leq C^\prime \big(H^{1+s} |f|_{H^s(\Omega)} + C_{\mu}^{1/2}(H,\ell)\ell^{d/2}(\sigma_\mu(H,\ell) + \mathsf{tol}  \ell^{d/2})\tnormf{f}{\Omega}\big),
	\end{split}
	\label{eq:errorest}
	\end{align}
	with $C_{\mu}$ from~\cref{a:Rieszbasis} and $|\cdot|_{H^s(\Omega)}$ denoting the $H^s$-seminorm.
\end{theorem}

From~\cref{eq:errorest}, we can identify the three ingredients that contribute to the overall RB-SLOD error, namely (i) the approximation of $f$ by the $L^2$-projection onto $\Pnull(\TH)$, (ii) the super-localization error, and (iii) the RB error.

\begin{proof}[Proof of \cref{th:error}]
	Let $\mu \in \ptr$ be arbitrary but fixed. 
	We first add and subtract the function $\bar u_\mu \coloneqq \mathcal A_\mu ^{-1} \PiH f \in \V_{H,\mu}$ and apply the triangle inequality to obtain
	\begin{equation*}
	\vnormof{u_\mu - u_{H,\ell,\mu}^\mathrm{rb}} \leq \vnormof{u_\mu - \bar u_\mu} + \vnormof{  u_{H,\ell,\mu}^\mathrm{rb}-\bar u_\mu}.
	\end{equation*}
	For the first term, one can show using the definition of $\bar u_\mu$,~\cref{eq:conditions4laxmilgram,eq:L2approx} that 
	\begin{align*}
	\alpha\vnormof{u_\mu-\bar u_\mu}&\leq   \sup_{v \in \V} \frac{A_\mu(u_\mu-\bar u_\mu,v)}{\vnormof{v}} = \sup_{v \in \V} \frac{\tspf{f-\PiH f}{v-\PiH v}{\Omega}}{\vnormof{v}}\\ 
	&\leq \pi^{-1}  H  \tnormf{f - \PiH f}{\Omega}.
	\end{align*}
	For the second term, we define
	\begin{equation*}
	\bar \varphi_{K,\ell,\mu}\coloneqq \mathcal A^{-1}_\mu g^\mathrm{rb}_{K,\ell,\mu}
	\end{equation*}
	and denote by $(c_{K})_{K\in \TH}$ the coefficients of the expansion of $\PiH f$ in terms of the right-hand sides $g_{K,\ell,\mu}^\mathrm{rb}$. Using the same coefficients, $\bar u_\mu$ can be expanded as 
	\begin{equation*}
	\bar u_\mu = \sum_{K \in \T_{H,\omega}} c_{K} \bar \varphi_{K,\ell,\mu}.
	\end{equation*}
	Abbreviating $e \coloneqq u_{H,\ell,\mu}^\mathrm{rb} - \bar u_{\mu}$, using~\cref{eq:conditions4laxmilgram} and the definition of the collocation solution~\cref{eq:urbslod}, we obtain 
	\begin{align}
	\label{eq:estlocerr}
	\alpha \vnormof{e}^2 \leq A_\mu(e,e) = \sum_{K \in \T_{H}}c_{K} A_\mu(\psi_{K,\ell,\mu}^\mathrm{rb}-\bar \varphi_{K,\ell,\mu},e).
	\end{align}
	
	We consider the summand corresponding to element $K \in \TH$ separately and denote the $\ell$-th order patch around $K$ by $\omega \coloneqq \mathsf N^\ell(K)$. We omit the fixed indices $K,\ell,$ and~$\mu$ of the functions $\bar \varphi_{K,\ell,\mu}$ and $\psi_{K,\ell,\mu}^\mathrm{rb}$ and denote by $\mathrm{tr}_\Sigma$ and $\mathrm{tr}^{-1}_\Sigma$ the trace and inverse trace operators with respect to $\Sigma = \partial \omega \backslash\partial \Omega$, respectively, see~\cref{eq:X,eq:conttrinv}. 
	Defining $\bar \psi \coloneqq \mathcal{A}_{\omega,\mu}^{-1} g^\mathrm{rb}$, using~\cref{eq:conormalder} and~$e-\operatorname{tr}^{-1}\operatorname{tr}e \in \V_\omega$, we obtain
	\begin{align}
	&A_\mu(\psi^\mathrm{rb}-\bar \varphi,e) =  A_\mu(\psi^\mathrm{rb},e)-\tspf{g^\mathrm{rb}}{e}{\omega} \notag\\
	&\qquad = A_{\omega,\mu}(\psi^\mathrm{rb},\mathrm{tr}_\Sigma^{-1}\mathrm{tr}_\Sigma e) -\tspf{g^\mathrm{rb}}{\mathrm{tr}^{-1}_\Sigma\mathrm{tr}_\Sigma e}{\omega}+ A_{\omega,\mu}(\psi^\mathrm{rb}-\bar \psi,e-\mathrm{tr}_\Sigma^{-1}\mathrm{tr}_\Sigma e)\notag\\
	&\qquad = \langle a_\mu\nabla \psi^\mathrm{rb}\cdot \nu,\mathrm{tr}_\Sigma e\rangle_{X^\prime\times X} + A_{\omega,\mu}(\psi^\mathrm{rb}-\bar\psi ,e-\mathrm{tr}_\Sigma^{-1}\mathrm{tr}_\Sigma e).\label{eq:lasteqproof}
	\end{align}
	For the bound of the first term in~\cref{eq:lasteqproof}, we denote by  $g \in \Pnull(\T_{H,\omega})$ the $L^2$-normalized SLOD right-hand side~\cref{eq:chocieg} and define 
	\begin{equation*}
	\tilde \psi^\mathrm{rb} \coloneqq \sum_{T \in \T_{H,\omega}} g|_T\chi_{T,\mu}^\mathrm{rb},\qquad\psi  \coloneqq\mathcal{A}_{\omega,\mu}^{-1} g
	\end{equation*}
	with $\chi_{T,\mu}^\mathrm{rb}$ defined in~\cref{eq:rbbestapprox}. Using~\cref{lemma:xprimenorm}, the definition of the conormal derivative~\cref{e:nd1}, the definitions of $g^\mathrm{rb}$ and $\sigma_\mu$ in~\cref{eq:chociegrb,eq:sigma}, respectively, as well as~\cref{eq:conditions4laxmilgram} and the continuity of $\mathrm{tr}^{-1}_\Sigma$ in~\cref{eq:conttrinv}, there holds
	\begin{align}
	\|a_\mu\nabla \psi^\mathrm{rb}\cdot \nu\|_{X^\prime} &= \vnormf{\mathcal R^\mathrm{rb}g^\mathrm{rb}}{\omega} \leq \vnormf{\mathcal R^\mathrm{rb} g\,}{\omega}\notag \\
	& = \sup_{\substack{w\in X\with \|w\|_X = 1}} |A_\mu(\tilde \psi^\mathrm{rb},\mathrm{tr}_\Sigma^{-1}w)-\tspf{g}{\mathrm{tr}_\Sigma^{-1} w}{\omega}|\notag \\
	&\leq \|a_\mu \nabla\psi \cdot \nu \|_{X^\prime} + \sup_{\substack{w\in X\with\|w\|_X = 1}} |A_\mu(\tilde \psi^\mathrm{rb}-\psi,\mathrm{tr}_\Sigma^{-1}w)|\notag \\
	&\leq \sigma_\mu(H,\ell) + \beta \,\vnormf{\tilde \psi^\mathrm{rb}-\psi}{\omega}.\label{eq:secondterm}
	\end{align}
	Using that the norm of $p$ in~\cref{eq:esterr} can be bounded as follows
	\begin{equation*}
	\vnormf{p}{\omega} = \|\one_T\|_{\V^\prime,\omega} \leq  \tnormf{\one_T}{\omega},
	\end{equation*} 
	we obtain for the second term in~\cref{eq:secondterm}
	\begin{align}
	\label{eq:estrberr}
	\begin{split}
	\vnormf{\tilde\psi^\mathrm{rb}-\psi}{\omega}
	&\leq \sum_{T \in \T_{H,\omega}} |g|_T|\vnormf{\chi_{T,\mu}^\mathrm{rb}-\chi_{T,\mu}}{\omega}\leq \alpha^{-1}\mathsf{tol}\sum_{T \in \T_{H,\omega}}g|_T\tnormf{\one_T}{\omega}\\&
	\leq C_\mathrm{p}\alpha^{-1}\mathsf{tol} \,\ell^{d/2},
	\end{split}
	\end{align}
	where $C_\mathrm{p}>0$ is some constant satisfying $\#\T_{H,\omega} \leq C_\mathrm{p}^2 \ell^d$. In estimate~\cref{eq:estrberr}, we employed~\cref{lemma:estimator}, abortion criterion~\cref{eq:esterr}, the discrete Cauchy--Schwarz inequality and the normalization  $\tnormf{g}{\omega} = 1$. 
	
	For the second term in~\cref{eq:lasteqproof}, we obtain, using the continuity of $\mathrm{tr}_\Sigma$ and $\mathrm{tr}^{-1}_\Sigma$ in~\cref{eq:conttr} and~\cref{eq:conttrinv}, respectively, that
	\begin{align*}
	A_{\omega,\mu}(\psi^\mathrm{rb} - \bar \psi,e-\operatorname{tr}^{-1}\operatorname{tr}e) &\leq \beta \vnormf{\psi^\mathrm{rb}-\bar \psi }{\omega}\vnormf{e-\operatorname{tr}^{-1}\operatorname{tr}e}{\omega}\leq 2\beta \vnormf{\psi^\mathrm{rb}-\bar \psi }{\omega}\vnormf{e}{\omega}.
	\end{align*}
	Here, the first term can be estimated similarly as~\cref{eq:estrberr} by
	\begin{align*}
	\vnormf{\psi^\mathrm{rb}-\bar \psi}{\omega} \leq C_\mathrm{p}\alpha^{-1} \mathsf{tol}\ell^{d/2} 
	\end{align*}
	using that $\tnormf{g^\mathrm{rb}}{\omega} = 1$. 
	
	Finally, combining the above estimates, we are able to continue estimate~\cref{eq:estlocerr}. Using~\cref{eq:L2stab} and~\cref{ass:stabg}, we obtain 
	\begin{align*}
	\alpha \vnormof{e}^2 &\leq \sum_{K \in \T_{H}}c_{K} A_\mu(\psi_{K,\ell,\mu}^\mathrm{rb}-\bar\varphi_{K,\ell,\mu},e)\\
	&\leq C C_\mathrm{ol} \ell^{d/2}\big(\sigma_\mu(H,\ell) +  C_\mathrm{p}\mathsf{tol}\beta\alpha^{-1} \ell^{d/2}\big) \sqrt{\sum_{K \in \T_{H}} c_{K}^2}\,\vnormof{e}\\
	&\leq C C_\mathrm{ol} C_{\mu}^{1/2}(H,\ell) \ell^{d/2}\big(\sigma_\mu(H,\ell) +  C_\mathrm{p}\mathsf{tol}\beta\alpha^{-1} \ell^{d/2}\big) \tnormf{f}{\Omega}\,\vnormof{e},
	\end{align*} 
	where $C_\mathrm{ol}>0$ reflects the overlap of the patches $\mathsf N^\ell(K)$ and $C>0$ is a generic constant independent of the discretization parameters. The assertion follows immediately.
\end{proof}

\begin{remark}[Choice of parameters]
	This remark investigates the choices of the oversampling parameter $\ell$ and the minimal number of elements $M$ in the sets $\prb$ required for achieving an error of order $H$ in~\cref{eq:errorest}. By the super-exponential decay of the quantity $\sigma_\mu$ as stated in~\cref{conj:sexpdec}, we need to choose $\ell$ of order $|\log H|^{(d-1)/d}$. 
	For the choice of $M$, let us recall from~\cref{rem:decdelta} that the greedy search error~\cref{eq:esterr} and, in particular, the achievable tolerance $\mathsf{tol}$, decay (sub-) exponentially in $M$. Denoting the exponent of~$M$ by $\gamma>0$, we need to choose~$M$ of order $|\log \ell H|^{1/\gamma}$.
\end{remark}

\section{Numerical experiments}
\label{sec:numexp}

For the numerical experiments, we consider the domain $\Omega = (0,1)^2$ equipped with a coarse Cartesian mesh $\TH$. Note that $H$ henceforth denotes the side-length of the mesh elements instead of their diameter. For any patch $\omega = \mathsf N^\ell(K)$, we discretize the corresponding local patch problems using the $\mathcal Q_1$-finite element method on the fine Cartesian meshes $\T_{h,\omega}$ obtained by uniform refinement of the coarse patch mesh~$\T_{H,\omega}$.
In the subsequent numerical experiments, we employ the version of the RB-SLOD as described in~\cref{sec:implementation}.

\subsection{Model diffusion problem}

The first numerical experiment is taken from~\cite[Sec.~4.1]{AbH15} and considers a pure diffusion problem with homogeneous Dirichlet boundary conditions. Its anisotropic multiscale diffusion tensor exhibits oscillations on various scales, the smallest being at the scale~$2^{-6}$ (see the reference for the exact definition of the coefficient). The problem's bilinear form admits an affine decomposition as in~\cref{eq:decA} with $Q = 4$ terms. The parameter space is the one-dimensional interval $\pset = [0,5]$ and, as training set $\ptr$, we use 100 equidistantly distributed points. 
For the fine-scale discretization of the patch problems, we use the mesh size $h = 2^{-8}$. 

\begin{figure}[ht]
	\begin{tabularx}{\textwidth}{@{}Y@{}}
		\includegraphics[width = .475\linewidth]{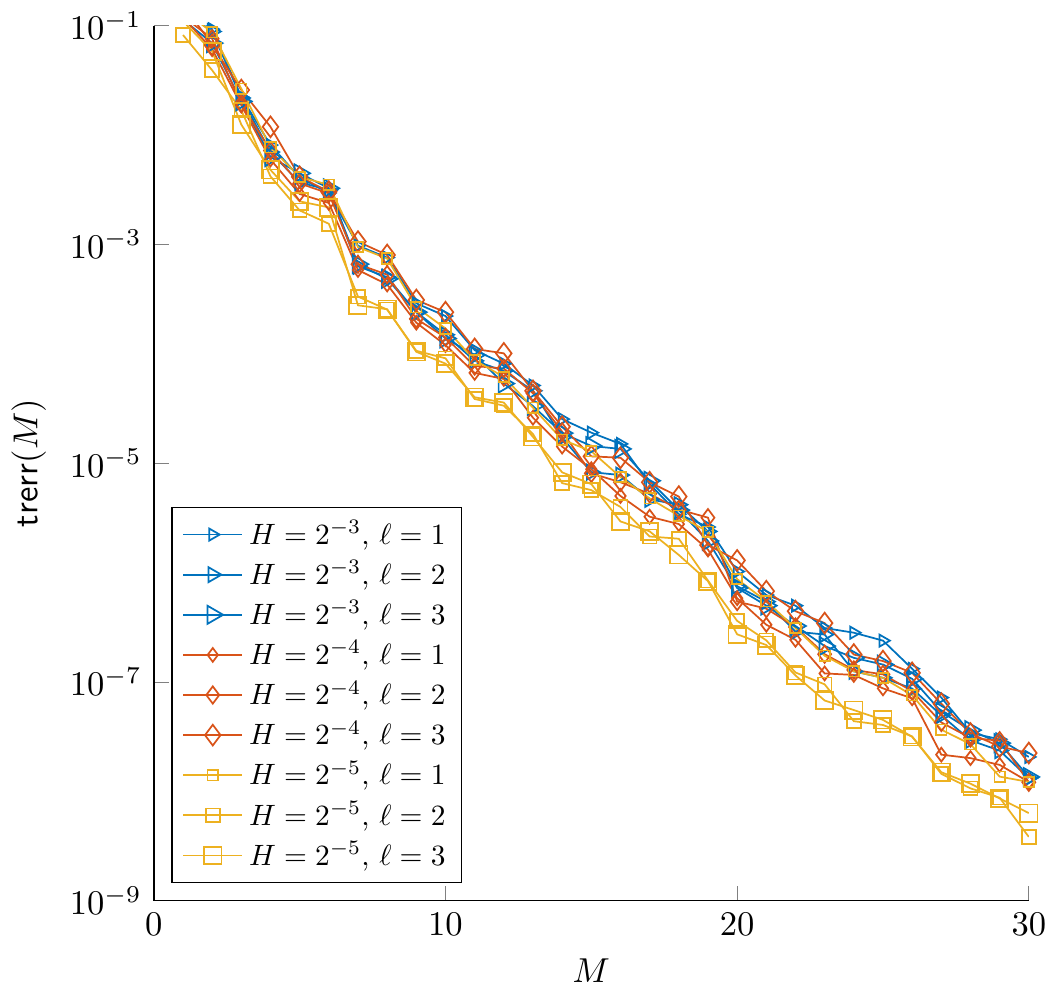}
	\end{tabularx}
	\caption{Decay of training error in greedy search algorithm for fixed patch $\omega$ and element $T$.}
	\label{fig:rbslodtraining}
\end{figure}

For demonstrating the decay of the training error in the greedy search algorithm, we fix a patch $\omega = \mathsf N^\ell(K)$ in the interior of $\Omega$ along with an element $T \subset \omega$. Recalling that the error estimator~$\Delta_\mu$ defined in~\cref{eq:Delta} depends on the number of elements $M$ in $\prb$, we define 
\begin{equation}
\label{eq:trerr}
\mathsf{trerr}(M) \coloneqq \frac{\max_{\mu \in \ptr} \Delta_\mu(M)}{\vnormf{p}{\omega}},
\end{equation}
which is the achievable tolerance of the greedy search algorithm ($p$ is defined in~\cref{eq:p}).
\Cref{fig:rbslodtraining} shows the decay of the training error for several pairs of $H$ and  $\ell$ as $M$ is increased.  In particular, we can observe a (sub-)exponential decay of the training error in $M$ confirming~\cref{rem:decdelta}. The discretization  parameters $H$ and $\ell$ have a small impact on the training error. Note that the qualitative decay behavior and the order of magnitude of the training errors are similar for all patches. 
Thus, given the discretization parameters $H$ and $\ell$ and a tolerance~$\mathsf{tol}$, one can read off the approximate size of the sets~$\prb$ necessary for reaching the desired tolerance from~\cref{fig:rbslodtraining}.

Next, we consider the constant right-hand side $f \equiv 1 \in \Pnull(\TH)$. Note that, for such right-hand sides, only the localization error and the RB error is present, i.e., the first term in~\cref{eq:errorest} vanishes. 
We compute the RB-SLOD solution for $\mu \approx 2.129$, which is not contained in the training set~$\ptr$. In order to highlighting that fine-scale discretizations of the patch problems~\cref{eq:chi,eq:defr,eq:sqk,eq:p} are involved, we henceforth denote the RB-SLOD solution by $u_{H,\ell,\mu,h}^\mathrm{rb}$. We compute the relative errors in the $\V$-norm with respect to the $\mathcal Q_1$-finite element reference solution $u_{\mu,h}$ on the global Cartesian mesh of mesh size $h = 2^{-8}$. As comparison, we also depicted the respective errors for the SLOD.
\begin{figure}[ht]
	\begin{tabularx}{\textwidth}{@{}YY@{}}
		\includegraphics[width=\linewidth]{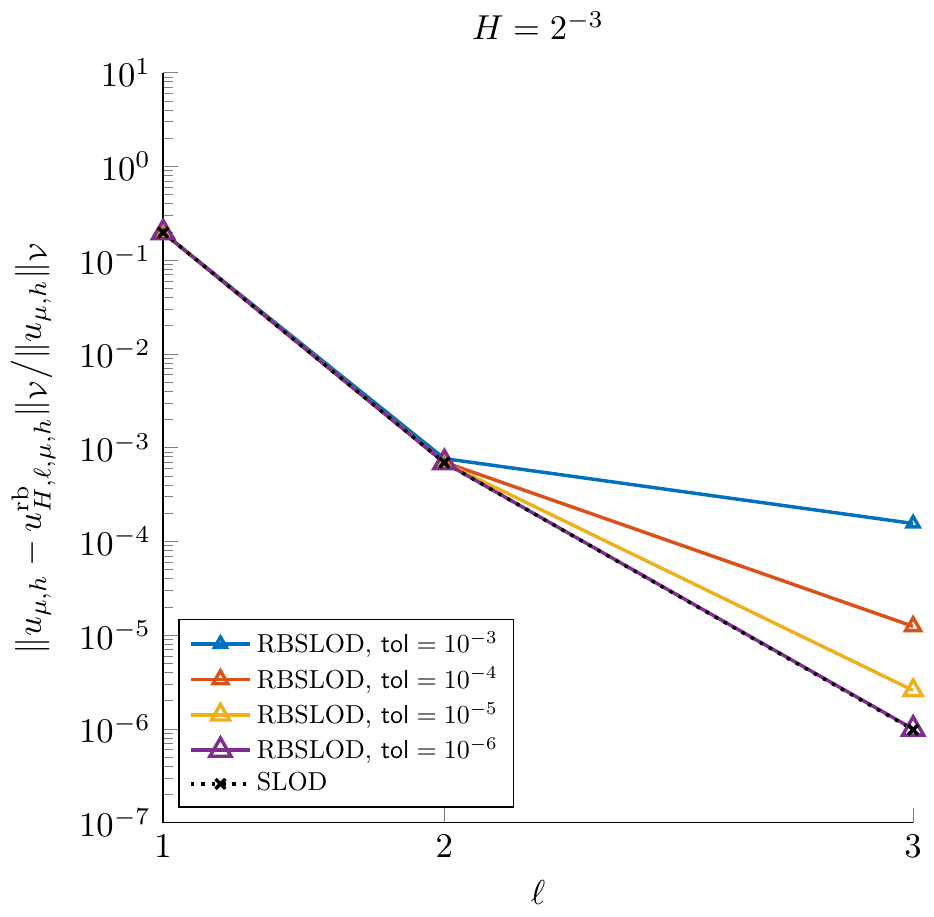}&
		\includegraphics[width=\linewidth]{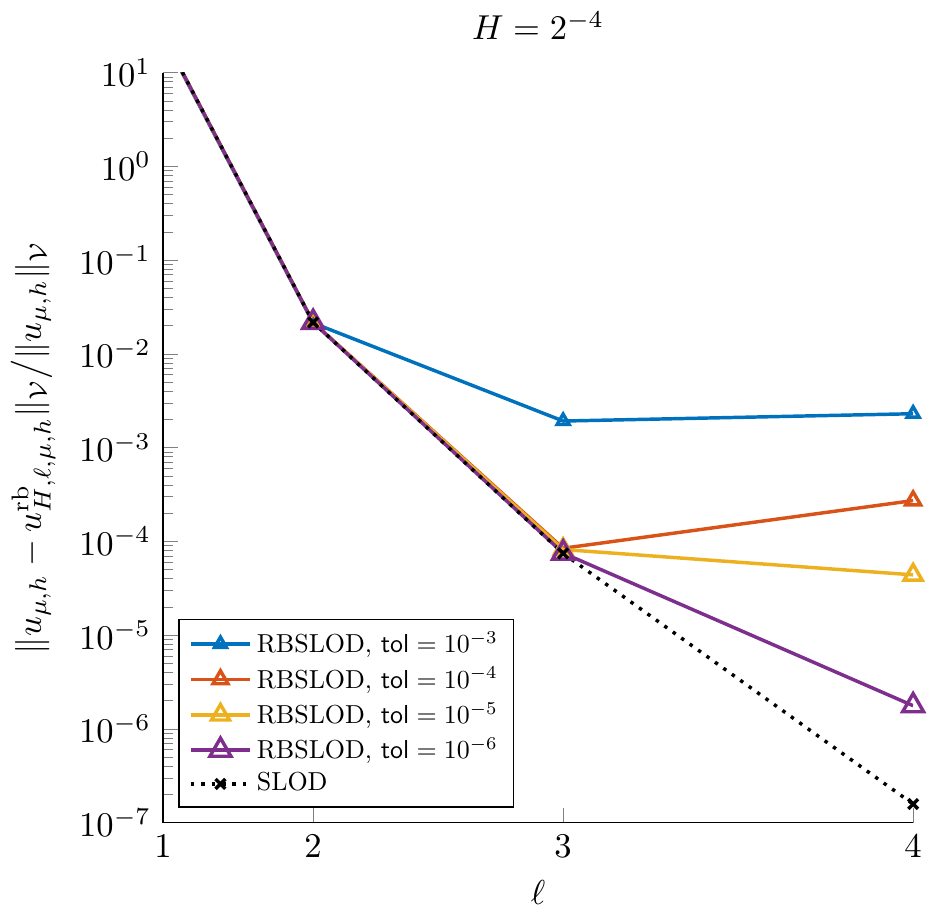}
	\end{tabularx}
	\caption{Relative errors of the RB-SLOD and SLOD for various tolerances of the greedy search versus the oversampling parameter~$\ell$ for $H = 2^{-3}$ (left) and $H = 2^{-4}$ (right). }
	\label{fig:rbsloderror}
\end{figure}
\Cref{fig:rbsloderror} shows that the error curve of the RB-SLOD approaches the error curve of the SLOD as $\mathsf{tol}$ is decreased. In the left plot, the error of the SLOD is reached for the value of the tolerance in greedy search $\mathsf{tol} = 10^{-6}$, whereas in the right plot, an even smaller tolerance would have been necessary.  
The error curve of the SLOD decays quadratic-exponentially (note that the quadratic scaling of the $x$-axis lets a quadratic-exponentially decaying curve appear linear), which is in line with the original results from~\cite{HaPe21b} and in agreement with \cref{conj:sexpdec}. 

\begin{figure}[ht]
	\begin{tabularx}{\textwidth}{@{}Y@{}}
		\includegraphics[width=.475\linewidth]{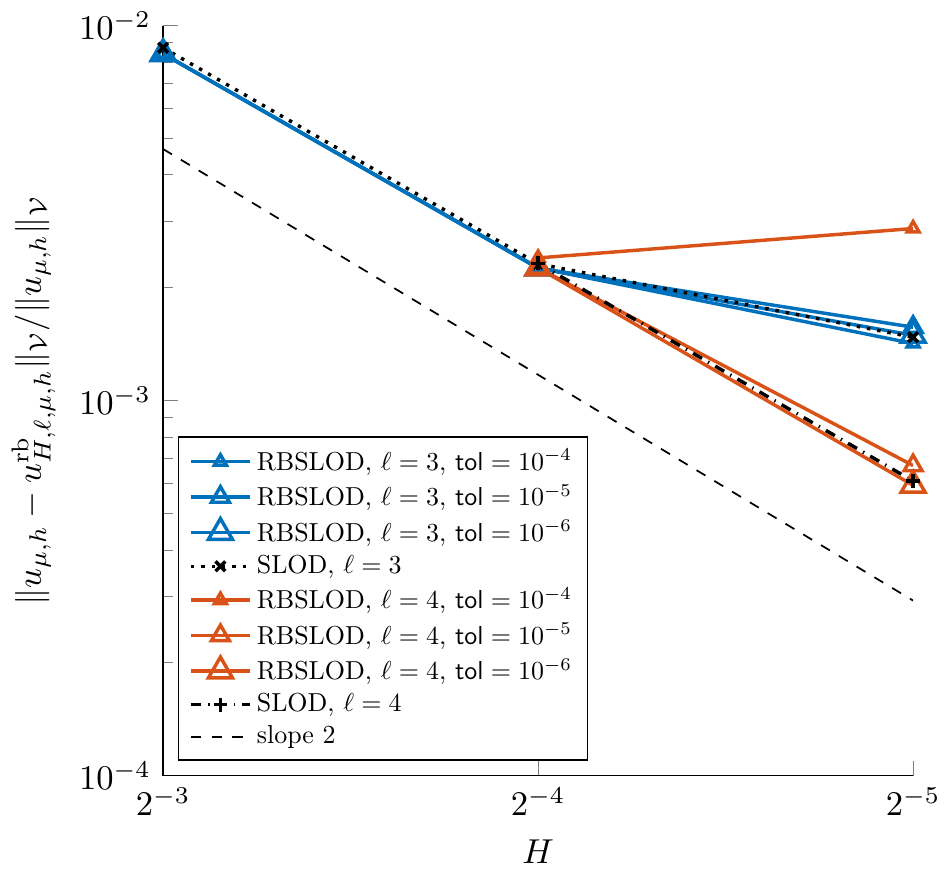}
	\end{tabularx}
	\caption{Relative $\V$-norm error for the RB-SLOD and SLOD in dependence of $H$ for several pairs of discretization parameters $\ell$ and $\mathsf{tol}$.}
	\label{fig:optimalorder}
\end{figure}

For confirming the method's optimal order of convergence in the mesh size $H$, as stated in~\cref{th:error}, we pick the right-hand side 
\begin{equation}
f(x_1,x_2) = \sin(x_1)\sin(x_2).\label{eq:rhssin}
\end{equation}
Note that for both right-hand sides, one can use the same reduced basis spaces and only the cheap coarse solve in~\cref{algo:coarsesolve} needs to be repeated, see~\cref{rem:rbinrhs}.
\Cref{fig:optimalorder} shows the relative errors of the RB-SLOD and the SLOD as a function of $H$. As reference, also a line of slope two is depicted which is the expected order of convergence, cf.~\cref{th:error} (recall that $f \in H^1(\Omega)$).
For sufficiently large oversampling parameters~$\ell$, the expected second order convergence is clearly visible. Again, it can be observed that the RB-SLOD error curve approaches the SLOD error curve as $\mathsf{tol}$ is decreased.

Note that a direct comparison of the proposed method to the RB-LOD from~\cite{AbH15} is difficult. More specifically, for the only considered right-hand side $f \equiv 1$ in~\cite{AbH15}, the SLOD is unbeatable, as only the localization error and the RB error are present. In comparison, for the RB-LOD, due to a different construction of the ansatz space, also the spatial approximation error is present. 
Hence, it seems reasonable to compare the magnitude of the errors of the RB-SLOD for the sinusoidal right-hand side~\cref{eq:rhssin} with the one of the RB-LOD for the constant right-hand side.  
It can be observed, that for the same magnitude of errors, the proposed method requires much smaller oversampling parameters than the RB-LOD. For example, given $H = 2^{-4}$,  both methods achieve a relative $\V$-norm error of order $10^{-2}$ with $\ell = 2$ for the RB-SLOD ($2.038\times 10^{-2}$) and $\ell = 4$ for the RB-LOD ($2.834 \times 19^{-2}$, cf.~\cite[Tbl.~1]{AbH15}). This allows for a sparser coarse system of equations and a more localized computation of the basis functions.

\subsection{Parameterized mass transfer with parameterized source}\label{sec:numexp2}

As second numerical experiment, we consider the example from~\cite[Ch.~8.4]{QMN15}. It is a reaction-convection-diffusion problem, where the magnitude of the diffusion coefficient, the direction of the advection, and a Gaussian source term are modeled using a five-dimensional parameter vector $\mu \in \mathbb R^5$ with components $\mu_1 \in [0.01,0.1]$, $\mu_2 \in [0,2\pi)$, $\mu_3, \mu_4,\in [0.25,0.75]$, and $ \mu_5  \in [0.1,0.25]$. More precisely, in strong form, the equation reads 

\begin{equation*}
\begin{aligned}
-\mu_1\Delta u - b_\mu\cdot \nabla u+u &= f_\mu\qquad  &&\text{in } \Omega \coloneqq (0,1)^2,\\
\mu_1 \nabla u \cdot \nu  &= 0 \qquad &&\text{on }\partial \Omega,
\end{aligned}
\end{equation*}
with
\begin{equation*}
b_\mu = (
\cos(\mu_2), \sin (\mu_2))^T,
\end{equation*}
and 
\begin{equation*}
f_\mu(x_1,x_2) = \exp\Big(-\frac{(x_1-\mu_3)^2 + (x_2-\mu_4)^2}{\mu_5^2}\Big).
\end{equation*}

This problem is challenging due to several reasons. First, the problem's nature is strongly dependent on the magnitude of diffusion $\mu_1$, i.e., the proposed reduced basis method must be able to deal with both the diffusion and convection-dominated regime at the same time. Second, the parameter space is relatively high-dimensional which makes RB methods typically quite expensive. Third, the chosen source term does not admit an affine decomposition. This prevents the straightforward application of standard RB approaches and, typically, further tools (like the empirical interpolation method~\cite{BMN04}) need to be employed, causing extra computational costs.

For such numerical examples, the RB-SLOD has the decisive advantage that it has a built-in RB approach in the right-hand side. More specifically, the right-hand side parameters are only important for the coarse solve in~\cref{algo:coarsesolve} and can be ignored in the actual reduced basis approach, see~\cref{sec:protapp,rem:rbinrhs}. In addition, as the considered coefficients are constant in space and in particular periodic with respect to~$\TH$, only $\mathcal O(\ell^d)$ patches need to be considered for the computations and one can deal with the remaining patches by translation, see~\cite[Rem.~B.1]{HaPe21b}. This drastically reduces the computational costs. 

\begin{figure}[ht]
	\begin{tabularx}{\textwidth}{@{}Y@{}}
		\includegraphics[width=.475\linewidth]{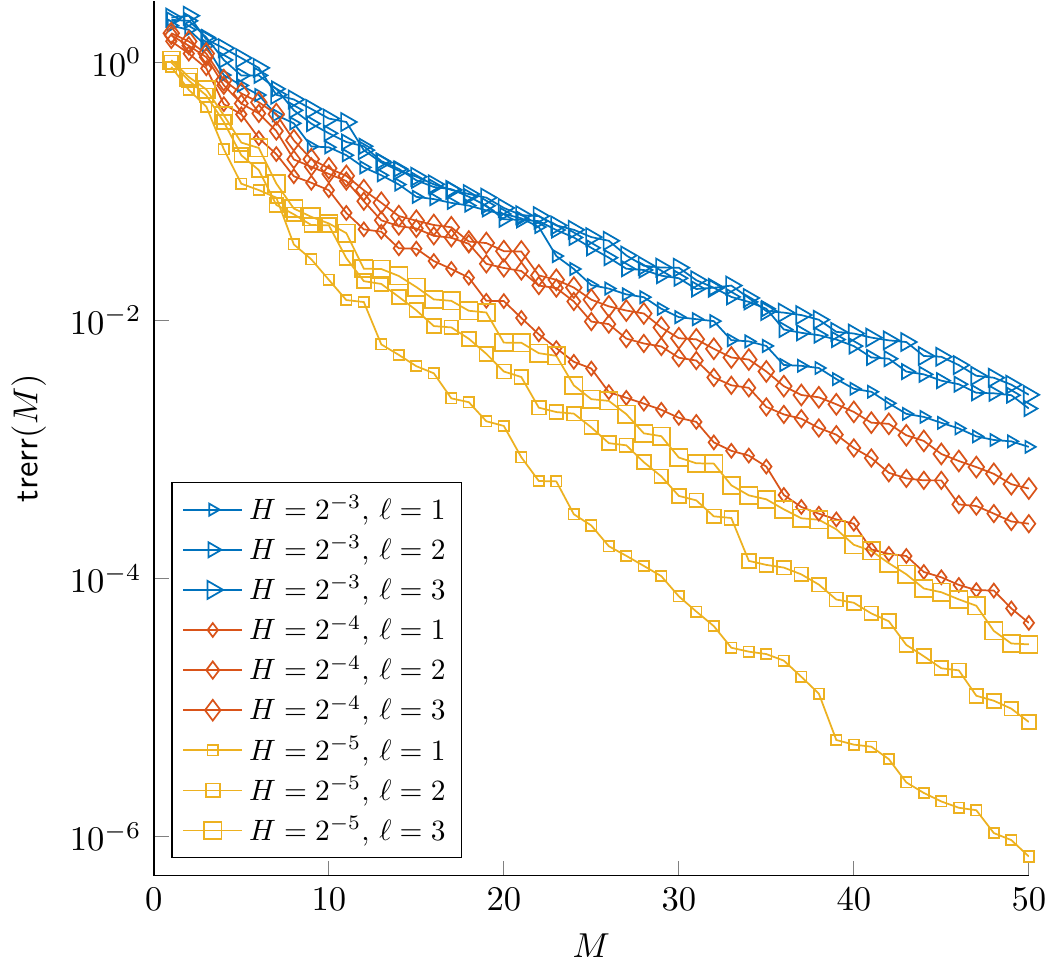}
	\end{tabularx}
	\caption{Decay of training error in greedy search algorithm for fixed patch $\omega$ and element $T$.}
	\label{fig:masstransfertrainingerr}
\end{figure}

In the implementation, we use a training set $\ptr$ of 400 points on a Cartesian grid of the two dimensional parameter set $\pset = [0.01,0.1]\times [0,2\pi]$, where we ignore the parameters stemming for the right-hand side. For the fine-scale discretization of the patch problems, we use the mesh size $h = 2^{-8}$. 
We again start by demonstrating the decay of the training error~\cref{eq:trerr} as the size $M$ of the set $\prb$ is increased, see~\cref{fig:masstransfertrainingerr}.
Compared to~\cref{fig:rbslodtraining}, the magnitude of training errors is significantly larger which is due to the more complicated nature of the underlying problem. Here, the magnitude of the error depends on the choice of $H$ and $\ell$ in the sense that for patches with a large diameter, also the training error is large.

\begin{figure}[h]
	\begin{tabularx}{\textwidth}{@{}YYYl@{}}
		$\mu_1 \approx 0.048,\; \mu_2 \approx 5.118$ & $\mu_1 \approx 0.090 ,\; \mu_2 \approx 5.391$ &  $\mu_1 \approx 0.027,\; \mu_2 \approx 2.046$\\[1ex]
		\includegraphics[trim=0 0 60 0, clip,width=\linewidth]{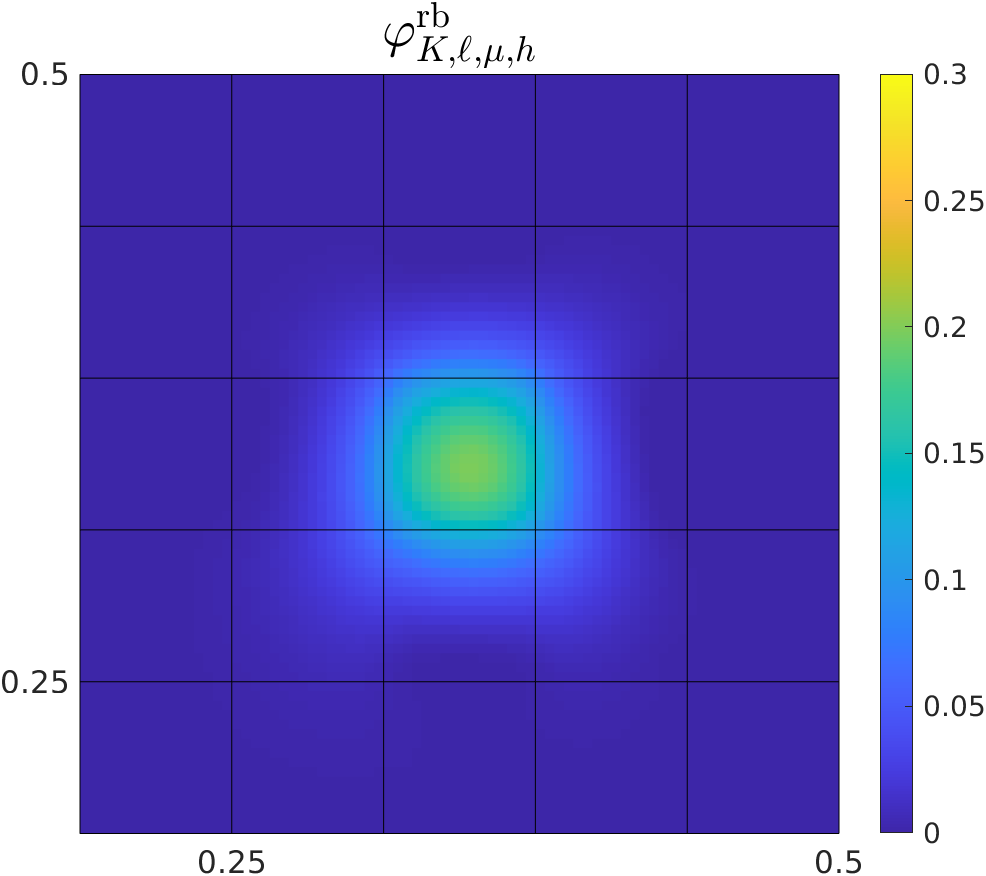}&
		\includegraphics[trim=0 0 60 0, clip,width=\linewidth]{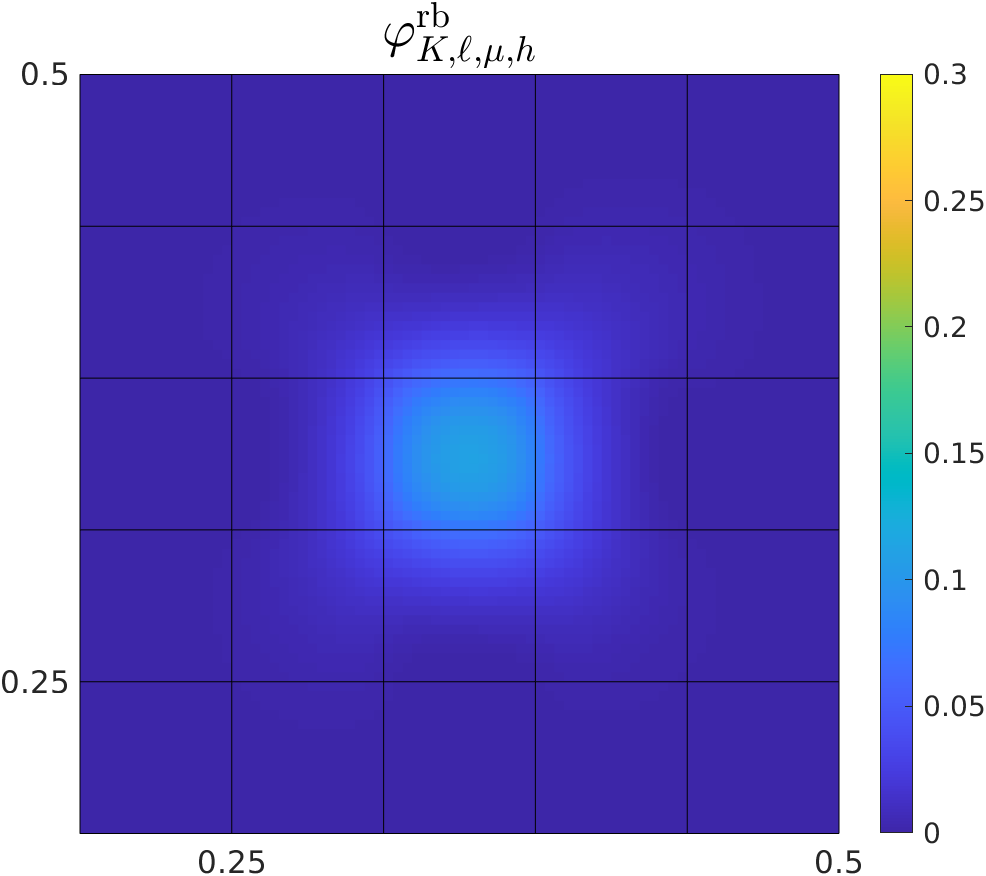}&
		\includegraphics[trim=0 0 60 0, clip,width=\linewidth]{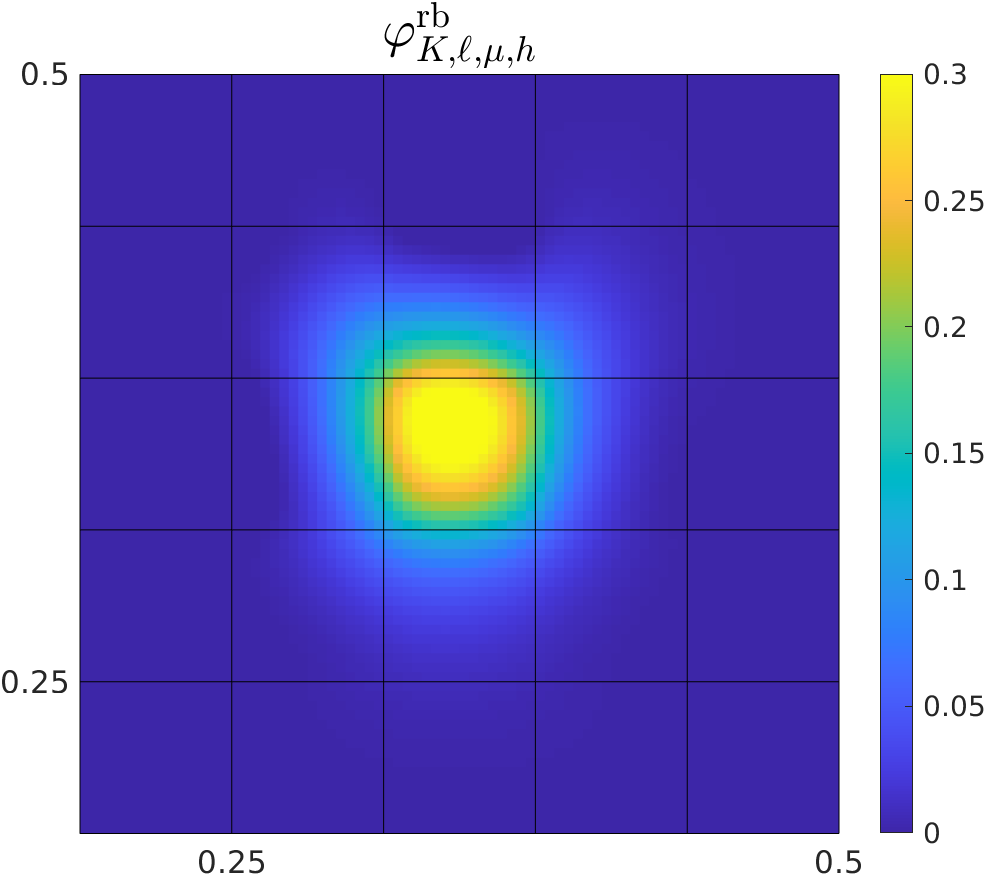} & \includegraphics[trim=422 0 0 0,clip,height=4.55cm]{figures/basisfunsamp1}
		\\[1ex]
		\includegraphics[trim=0 0 47 0, clip,width=\linewidth]{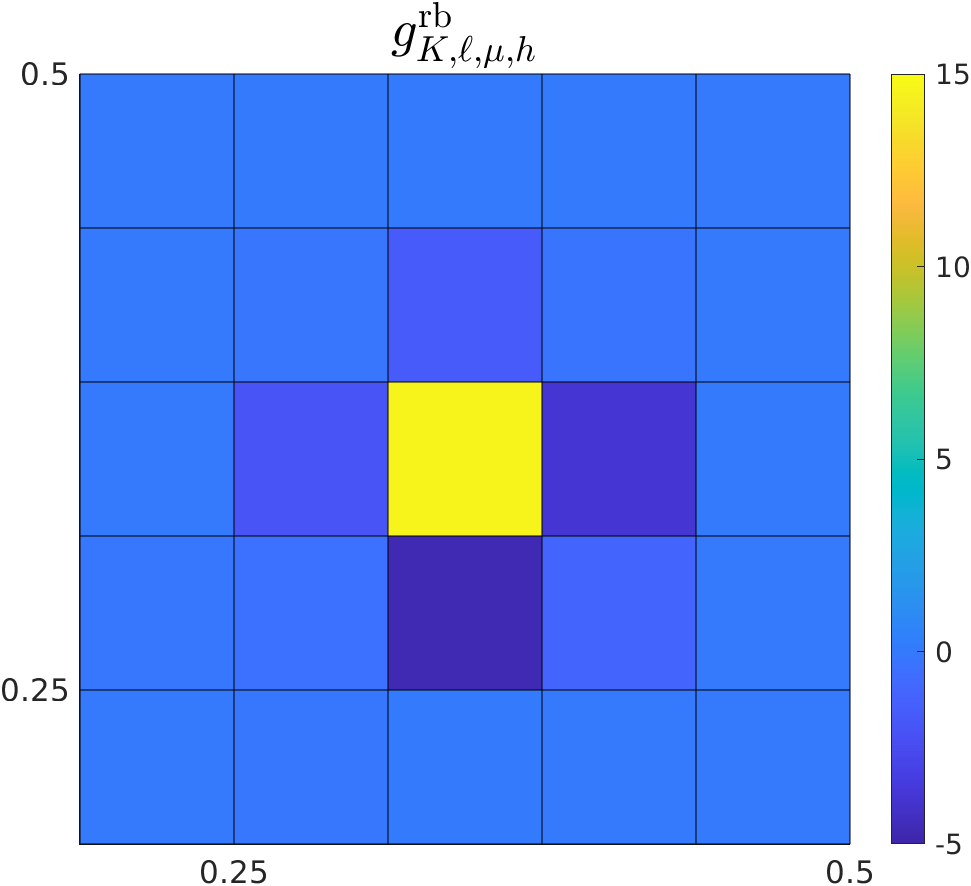}&
		\includegraphics[trim=0 0 47 0, clip,width=\linewidth]{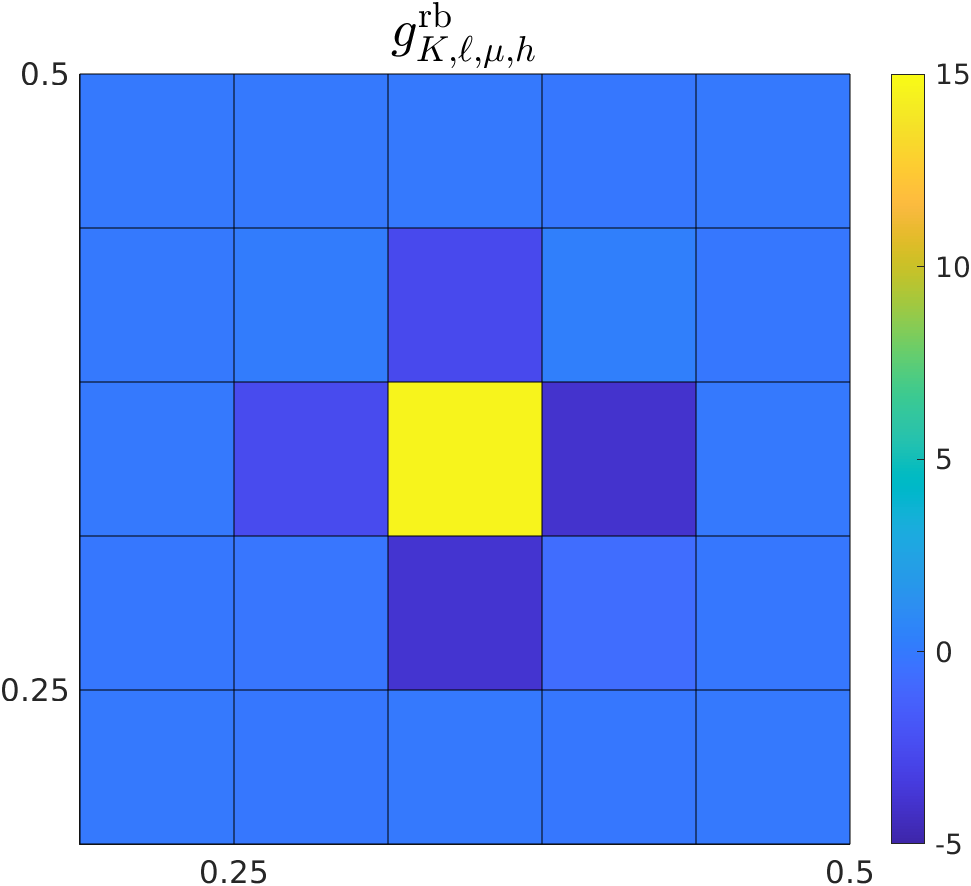}&
		\includegraphics[trim=0 0 47 0, clip,width=\linewidth]{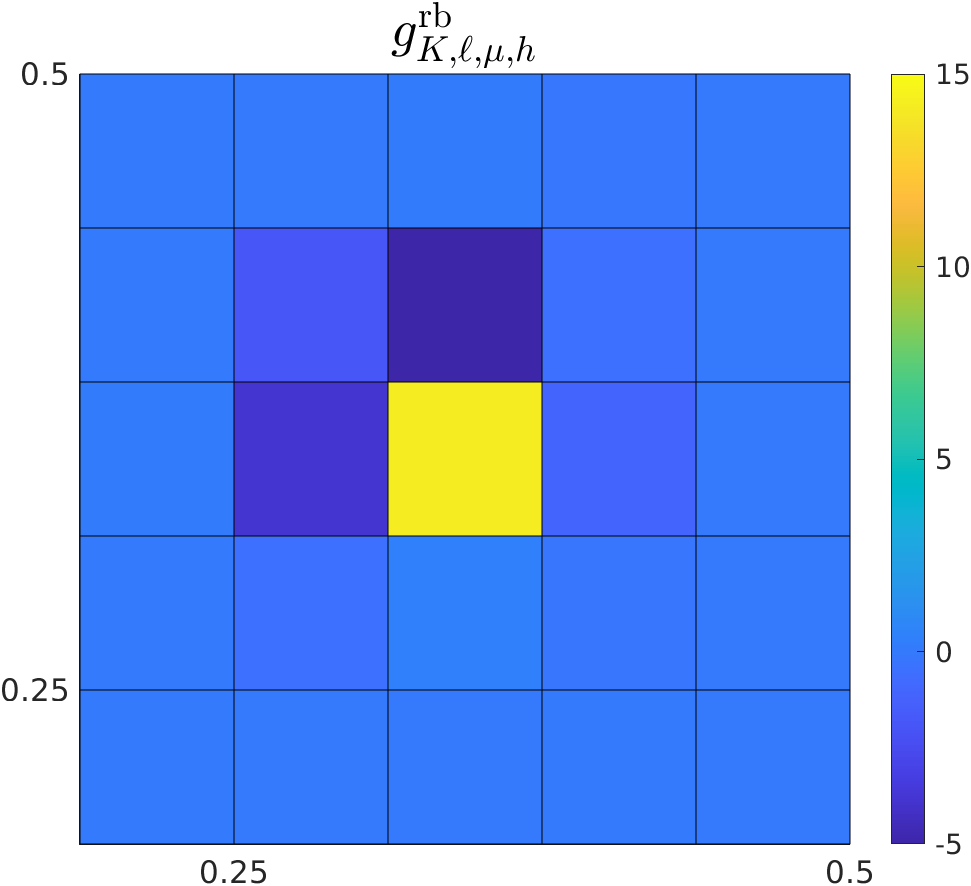} & 
		\includegraphics[trim=430 0 0 0, clip,height=4.55cm]{figures/rhssamp1}
	\end{tabularx}
	\caption{Depiction of RB-SLOD basis functions (top) and the corresponding right-hand sides (bottom) for three different parameter pairs (left to right).}\label{fig:slodbasis}
\end{figure}

\begin{figure}[h]
	\begin{tabularx}{\textwidth}{@{}YYYr@{}}
		$\mu_1 \approx 0.048,\, \mu_2 \approx 5.118,\linebreak \mu_3 \approx 0.512,\, \mu_4 \approx 0.703,\linebreak \mu_5 \approx 0.140$ &
		$\mu_1 \approx 0.090 ,\, \mu_2 \approx 5.391,\linebreak \mu_3 \approx 0.293,\, \mu_4 \approx 0.286,\linebreak\mu_5 \approx 0.130$ &
		$\mu_1 \approx 0.027,\, \mu_2 \approx 2.046,\linebreak \mu_3 \approx  0.608,\,\mu_4 \approx    0.703,\linebreak \mu_5 \approx 0.136$\\[1ex]
		\includegraphics[trim=0 0 60 0, clip,width=\linewidth]{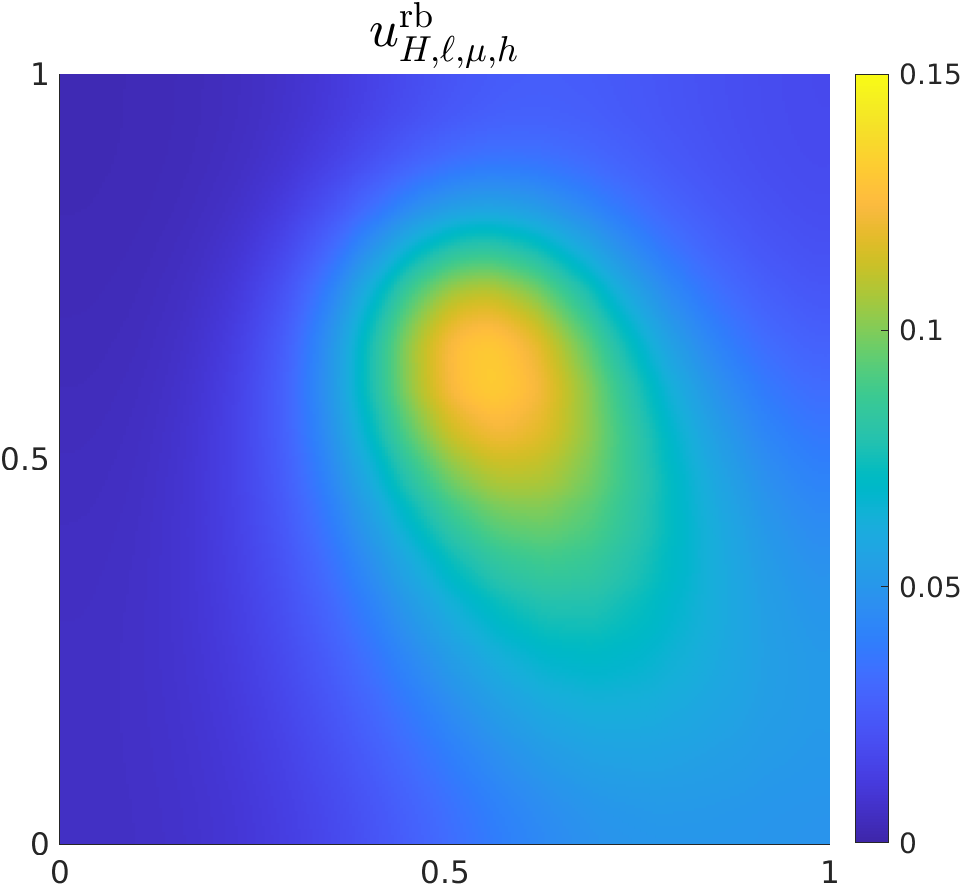}&
		\includegraphics[trim=0 0 60 0, clip,width=\linewidth]{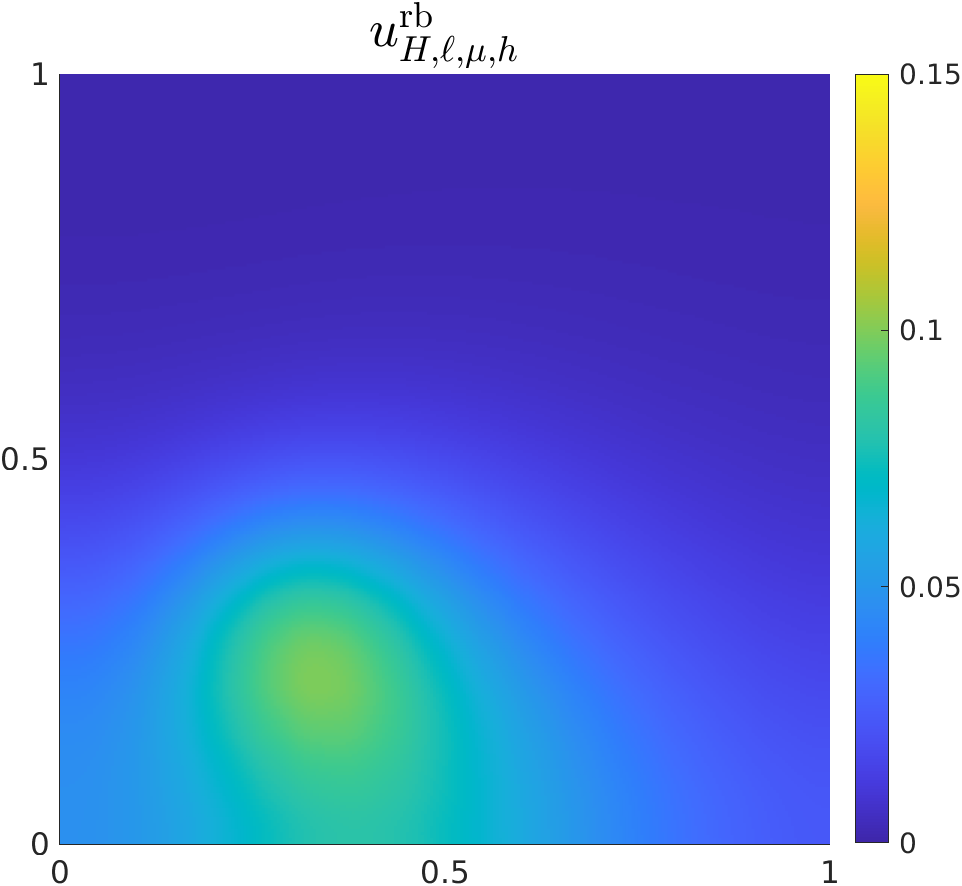}&
		\includegraphics[trim=0 0 60 0, clip,width=\linewidth]{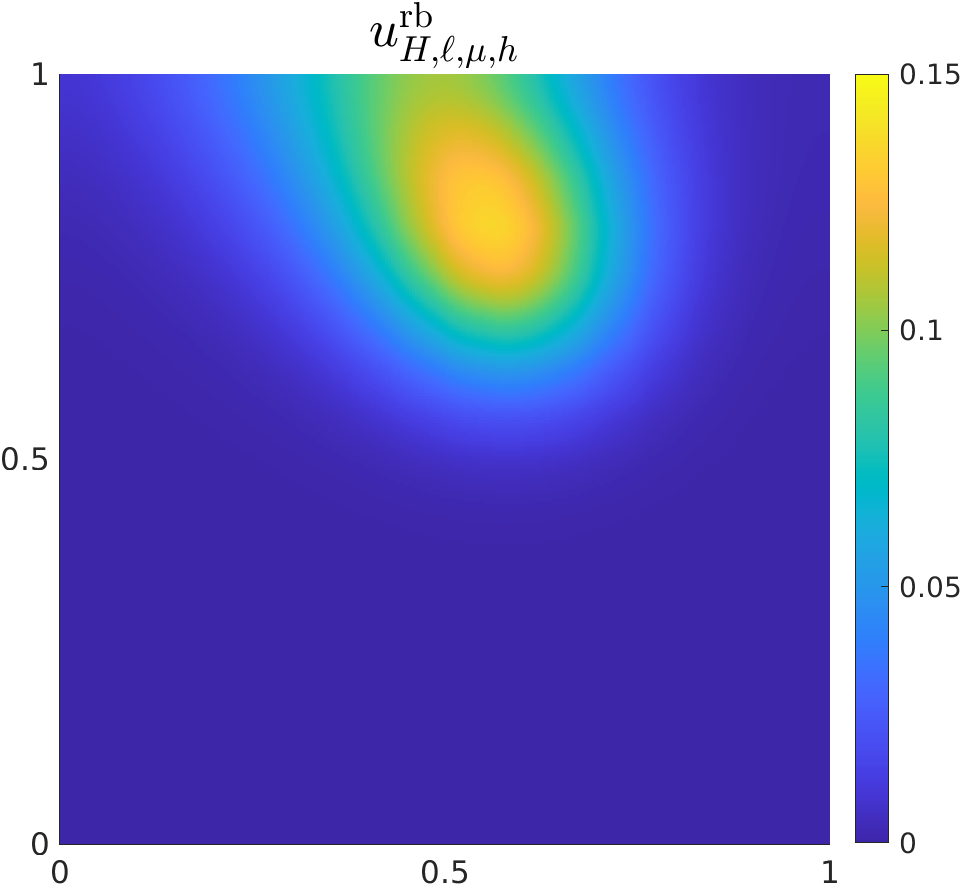}&
		\includegraphics[trim=408 0 0 0,clip,height=4.7cm]{figures/solsamp1}
		\\[1ex]
		\includegraphics[trim=0 0 60 0, clip,width=\linewidth]{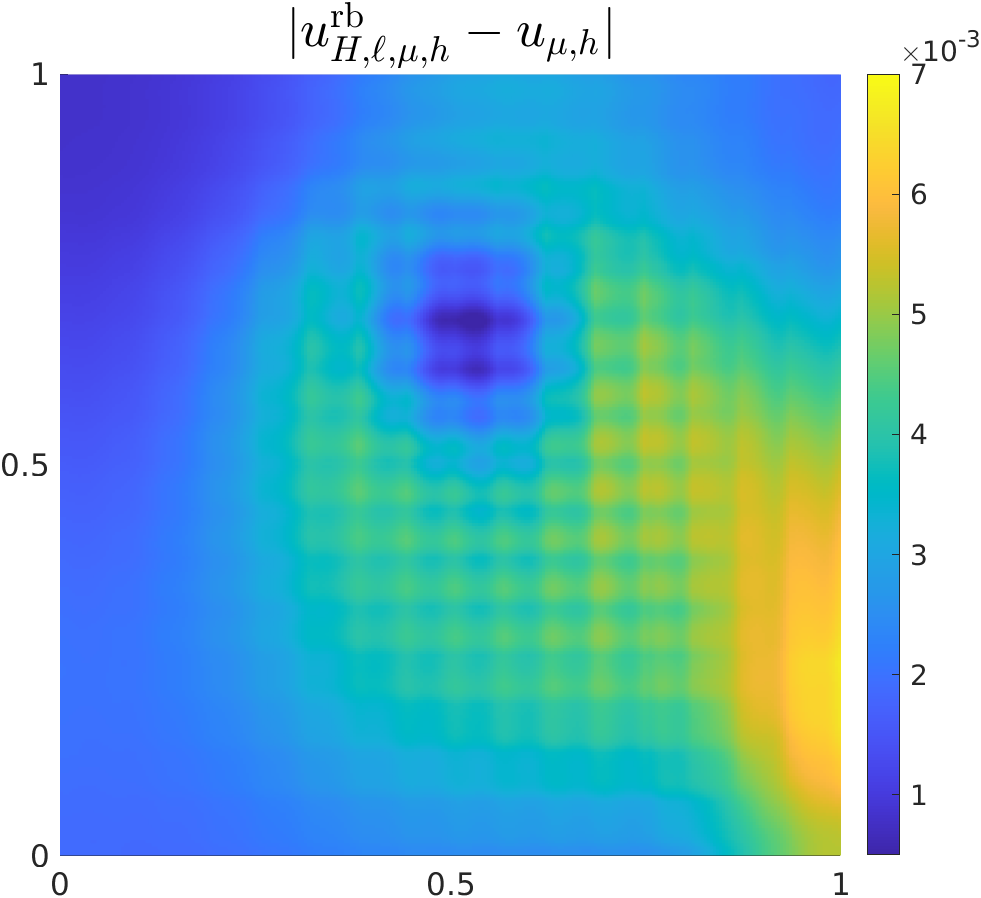}&
		\includegraphics[trim=0 0 60 0, clip,width=\linewidth]{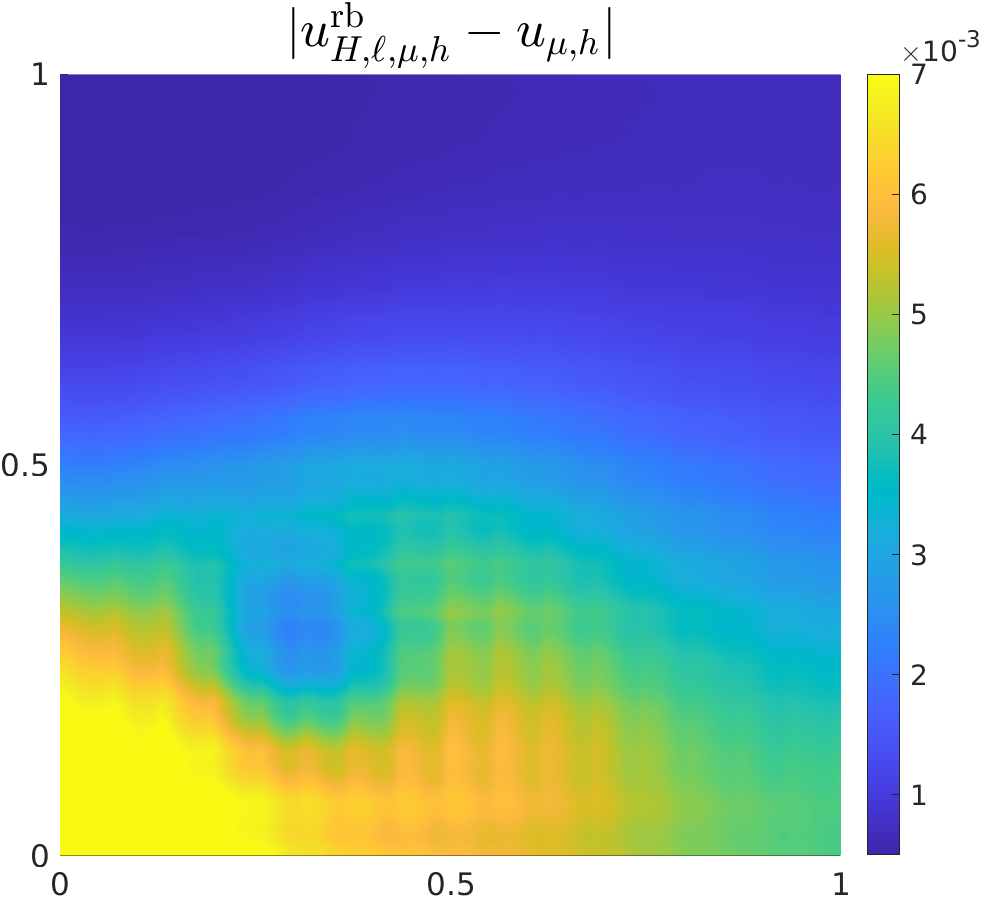}&
		\includegraphics[trim=0 0 60 0, clip,width=\linewidth]{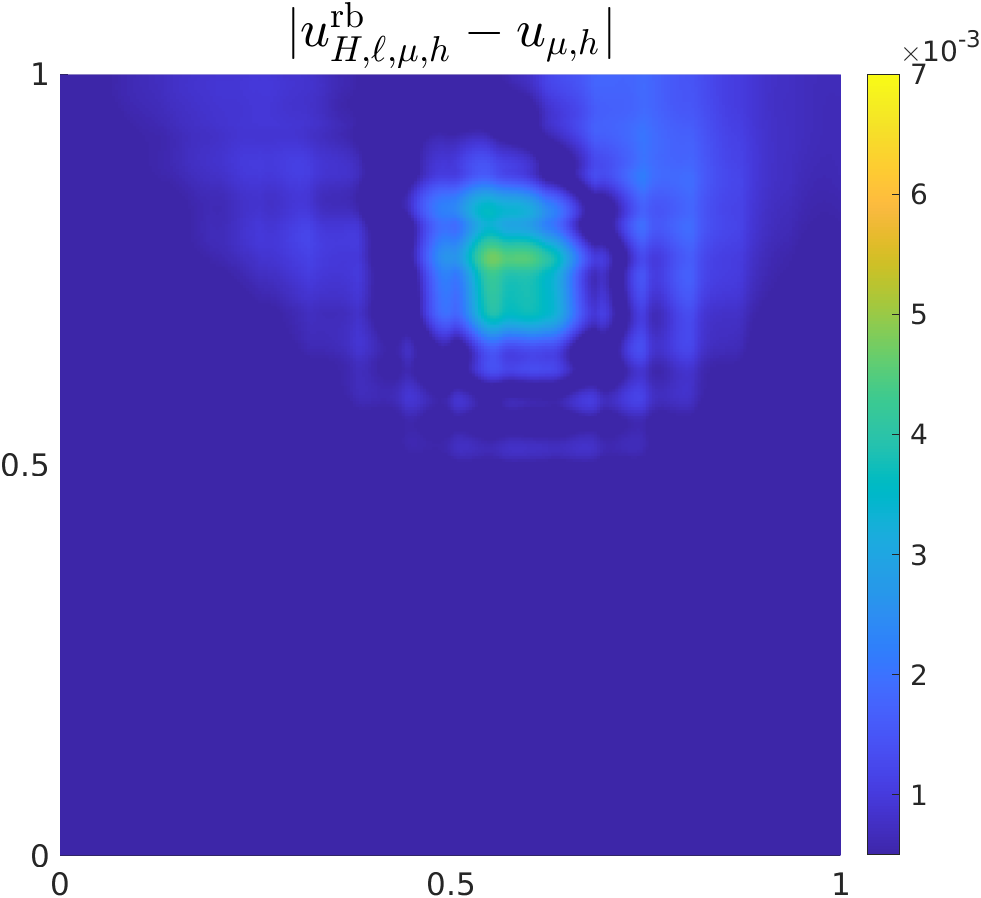} &
		\includegraphics[trim=415 0 0 0, clip,height=4.65cm]{figures/errsamp1}
	\end{tabularx}
	\caption{Depiction of RB-SLOD solutions (top) and the absolute value of the error (bottom) for three different parameter pairs (left to right).}
	\label{fig:slodsol}
\end{figure}

Next, we consider a coarse mesh with mesh size $H = 2^{-4}$ and the oversampling parameter $\ell = 2$. For this choice of discretization parameters, the largest patches make up at most 10\% of the whole domain's volume and the maximal number of elements in the sets $\prb$ is~34.
In~\cref{fig:slodbasis}, we depict one RB-SLOD basis function and the corresponding right-hand side for three different parameter values of $\mu_1$ and $\mu_2$. Note that the basis functions are independent of $\mu_3,$ $\mu_4,$ and $\mu_5$.
\Cref{fig:slodsol} depicts the RB-SLOD solution and the absolute value of the error with respect to the reference solution $u_{h,\mu}$ for the same parameter pairs as in~\cref{fig:slodbasis} and some random choices of the parameters $\mu_3,$ $\mu_4,$ and $\mu_5$. Note that all plots in the same row have the same color scale.
From left to right, the relative errors in the $\V$-norm are $5.62 \times 10^{-2}$, $4.53 \times 10^{-2}$, and $6.29 \times 10^{-2}$.
For this numerical experiment, a direct comparison to~\cite[Ch.~8.4]{QMN15} is again difficult. Therein, in the numerical experiments, the parameters $\mu_1$ and $\mu_5$ are fixed which simplifies the resulting problem. In particular, the nature of the problem does no longer depend on the choice of parameters (diffusion $\mu_1$ is fixed). In addition, the proposed method is able to avoid the extra difficulties with non-affine right-hand sides $f$ which classical reduced basis approaches have.

\section{Conclusion and outlook}
\label{sec:conclusions}

We have presented a novel method for the efficient approximate solution of parametric and strongly heterogeneous reaction-convection-diffusion problem. It combines the main principles of RB and SLOD: the RB is used for accelerating the typically costly SLOD basis computation, while the SLOD is employed for an efficient compression of the problem's solution operator requiring coarse solves only. 
Numerical experiments have confirmed the superiority of the proposed method, when compared to other state-of-the art techniques, like the RB-LOD. The super-localization properties of the SLOD allow one to perform the local computations on significantly smaller patches (for the same accuracy) reducing the computational costs considerably. Further, the more local support of the basis functions implies enhanced sparsity properties of the coarse system matrix. Given a value of the parameter vector, the method outputs a corresponding compressed solution operator which can be used for the efficient treatment of multiple, possibly non-affine, right-hand sides, only requiring one coarse solve per right-hand side. 

A natural follow up of the present work concerns the study of wave propagation and scattering problems in highly heterogeneous media. For such target, the SLOD has been developed in~\cite{Freese-Hauck-Peterseim} and model order reduction techniques for the frequency response map are studied in~\cite{Bonizzoni-Nobile-Perugia,Bonizzoni-Nobile-Perugia-Pradovera-a,Bonizzoni-Nobile-Perugia-Pradovera-b,Bonizzoni-Pradovera,Bonizzoni-Pradovera-Ruggeri}. 
Moreover, the possible application of the SLOD to improve numerical stochastic homogenization methods~\cite{GaP19,Fischer-Gallistl-Peterseim,FeP20} will be object of future investigation.

\section*{Acknowledgments}
We would like to thank Christoph Zimmer for inspiring discussion on the localization problem and especially its implementation.

\bibliographystyle{alpha}
\bibliography{bib}
\end{document}